\documentclass[a4paper]{article}

\author{Joaquim Ro\'e i Vellv\'e}

\usepackage[centertags]{amsmath}
\usepackage{amscd}
\usepackage{amsthm}
\usepackage{amssymb}
\usepackage{xypic}
\usepackage[dvips]{graphics}

\theoremstyle{definition}
\newtheorem{Def}{Definition}[section]
\theoremstyle{plain}
\newtheorem{Lem}[Def]{Lemma}
\newtheorem{Cor}[Def]{Corollary}
\newtheorem{Pro}[Def]{Proposition}

\theoremstyle{remark}
\newtheorem{Rem}[Def]{Remark}
\newtheorem{Exa}[Def]{Example}

\newlength{\taulength}
\makeatletter
\@addtoreset{table}{section}
\makeatother

\newcommand{\BF}{\operatorname{BF}}
\newcommand{\Bl}{\operatorname{Bl}}

\newcommand{\Spec}{\operatorname{Spec}}
\newcommand{\Pic}{\operatorname{Pic}}

\newcommand{\rts}{\operatorname{roots}}
\newcommand{\frs}{\operatorname{frees}}
\newcommand{\mult}{\operatorname{mult}}

\newcommand{\Isom}{{\bf Isom }}
\newcommand{\Cl}{{\mathit{Cl}}}
\newcommand{\Eff}{{\mathit{Eff}}}
\newcommand{\B}{{\bf B }}
\newcommand{\F}{{\bf F }}
\newcommand{\D}{{\bf D}}
\newcommand{\Du}{{\mathbb{D}}}

\newcommand{\Sch}{{{\cal S}ch}}
\newcommand{\Sets}{{{\cal S}ets}}

\newcommand{\Z}{\mathbb{Z}}

\newcommand{\m}{{\bf m}}
\renewcommand{\P}{\mathbb{P}}

\renewcommand{\O}{{\cal O}}
\renewcommand{\L}{{\cal L}}

\begin{document}

\title{Varieties of clusters and Enriques diagrams}
\maketitle

\begin{abstract}
Given a surface $S$ and an integer $r \ge 1$, there
is a variety $X_{r-1}$ parametrizing all clusters of
$r$ proper and infinitely near points of $S$
(see \cite{Kle77}). 
We study the geometry of the varieties $X_r$,
showing that for every Enriques diagram $\D$
of $r$ vertices the subset $\Cl(\D) \subset X_{r-1}$
of the clusters with
Enriques diagram $\D$ is locally closed.
We study also the relative positions of the
subvarieties $\Cl(\D)$, showing that they
do not form a stratification and giving
criteria for adjacencies between them.
\end{abstract}

\section*{Introduction}

In the 1970's S. Kleiman introduced
the method of iteration in multiple-point theory
for the study of maps $f:X \longrightarrow Y$ 
between algebraic varieties.
The basics of this method are set in \cite{Kle81a} (see 
also \cite{Kle77} and \cite{Kle82}), where the derived 
varieties $X_r$ and maps $X_r \longrightarrow X_{r-1}$, 
$r \in \Z_{\ge 0}$ are defined using blowing-ups
in a recursive way.
The varieties $X_r$ are thought of as parametrizing sets of
points of $X$, possibly infinitely near,
whose image under $f$ is a single point.

If $Y=\Spec k$ and $f:X \longrightarrow Y$ is
a smooth morphism, then the iterated blowing-ups 
parametrize ordered clusters of points of $X$. This was
implicit in Kleiman's work, and proved
or used afterwards by others in various settings (see
\cite{Har88}, \cite{Ran85}, \cite{DO88}). 

Similarly to the varieties $X_r$, other
varieties can be defined to parametrize 
unibranched  clusters
\cite{Ran85} or ordered unions of unibranched clusters
\cite{Col87}. In these papers Z. Ran and S. Colley
respectively compute fundamental classes of the
involved varieties as well, and use them to
enumerate \emph{united} and \emph{stationary}
multiple points of maps respectively. Similar work
for cases with three points was done by
J. Roberts in \cite{Rob80}.
The paper \cite{Ran85} contains also information on
the relation between varieties of clusters and the Hilbert scheme
of zero-dimensional subschemes of $S$; other works dealing
with this relation are \cite{Bel92}, 
\cite{Kee93}, and \cite{Evain}.

Here we study the geometry of the varieties 
of clusters of points on a surface, and some subvarieties
in them corresponding to clusters with
given Enriques diagrams (the combinatorial data
of proximities between points of a cluster $K$ are
encoded by means of \emph{Enriques diagrams}; two clusters have the same
Enriques diagram if and only if their points satisfy
the same proximity relations). The varieties
for unibranched clusters used by Z. Ran
and S. Colley arise as particular cases
of such subvarieties.

The paper is organized as follows. In section
\ref{pre} we review definitions and basic facts
on clusters, Enriques diagrams and varieties of clusters. 
In section \ref{fix} we study the
subvarieties of $X_r$ determined by Enriques
diagrams; given an ordered Enriques diagram $\D$ of $r$ vertices, we denote
$\Cl(\D) \subset X_{r-1}$ the set of all clusters whose
Enriques diagram is $\D$. In Propositions \ref{efclo} and \ref{loclo}
we prove that $\Cl(\D)$ is locally closed in $X_{r-1}$.
Next we consider the relative positions of the
varieties $\Cl(\D)$; in particular, given two
ordered Enriques diagrams $\D$ and $\D'$ we are
interested in giving combinatorial criteria
to decide whether $\Cl(\D') \subset \overline{\Cl(\D)}$
or not. Lemma \ref{noiguals} gives a necessary condition,
and in the remains of section \ref{fix} we
study under which conditions it is also a sufficient
condition.
In section \ref{nostrata} we show that the
varieties $\Cl(\D)$ do not form
a stratification, i.e. it is not true that
$\Cl(\D')\cap \overline{\Cl(\D)} \ne \emptyset$
implies $\Cl(\D') \subset \overline {\Cl(\D)}$.
In section \ref{functorial} we study families of 
blown-up surfaces, using the universal property of $X_{r-1}$ 
of \cite{Har88} together with results of the previous sections. 
Concretely, in Proposition \ref{diagsfam2} we prove that for
every proper and smooth morphism
$\psi: X \rightarrow T$, whose
fibres are blowing-ups of clusters
of $r$ points of $S$, and every $t \in T$, there exist 
an \'etale neighborhood $V$ of $t$ and an open dense subset 
$U$ of $V$ where the Enriques diagram of the
blown-up cluster is constant, and a relation
between this constant Enriques diagram and
that of the cluster blown up in $X_t$. This
gives a sort of ``semicontinuity'' of the
Enriques diagram.
Finally, in section \ref{puntfix} we study
subvarieties of $X_{r-1}$, defined by fixing
the position of some points of the cluster
in different ways.

Apart from their own interest, results of
this paper can be applied to the study of linear
systems (see \cite{Roe01}, \cite{Roe?2}, 
\cite{Har01}, \cite{Roe?3}, \cite{HR??}),
of maps from $X_r$ to the Hilbert scheme of fat points
of $S$ (see \cite{Roe?2}), or adjacencies of
equisingularity types (\cite{AR??}). We also
expect that they will have applications
to enumerative geometry, extending results of
\cite{Ran85}, \cite{Col87}, \cite{KP99}.

\section{Preliminaries}
\label{pre}
\label{vc}
\label{itera}

  Let $k$ be an algebraically closed field
and let $S$ be a smooth irreducible projective 
surface, defined over $k$. We begin by reviewing
some well-known facts dealing with clusters of
points of $S$, referring to the book \cite{Cas00}
for proofs and more on the subject.

Given a point $p$ of $S$, we denote the
\emph{blowing-up} of $p$ on $S$ by
$\pi_p:S_p \longrightarrow S$. The exceptional divisor $E_p$
of $\pi_p$ is called the \emph{first
(infinitesimal) neighbourhood} of $p$ (on $S$) and its points
are the points in the first neighbourhood of
$p$. If $i>0$, one defines by induction the \emph{points in
the $i$-th neighbourhood} of $p$ as the
points in the first neighbourhood of some
point in the $(i-1)$-th neighbourhood.
The points which are in the $i$-th neighbourhood of $p$
for some $i>0$, are also called points \emph{infinitely
near} to $p$. Sometimes the points in $S$ will be called
\emph{proper points} in order to distinguish them from
the infinitely near ones.

Let $p, q$ be two points in $S$, proper or infinitely near.
We will say that $p$ precedes $q$, $p < q$, if and only if
$q$ is infinitely near to $p$. We will write $p \le q$ if
$q$ is equal or infinitely near to $p$. The relation 
$\le$ is a partial ordering of the set of all points, proper
and infinitely near, called the \emph{natural
ordering}. A \emph{cluster} of points of a surface $S$ is a finite
set of points $K$ proper or infinitely near in $S$, such
that, for each point $p\in K$, $K$ contains all points
to which $p$ is infinitely near. 
We denote $S_K$ the surface obtained by blowing up all 
points in $K$, and $\pi_K:S_K \longrightarrow S$ the 
composition of the blowing-ups.
An \emph{admissible ordering} of a cluster $K$ is
a total ordering of the points of $K$, refining
the natural ordering, and an \emph{ordered cluster} is a cluster
together with an admissible ordering. Equivalently,
an ordered cluster may be defined as a sequence of
points $(p_1, p_2, \dots, p_r)$ such that $p_1$
is a proper point of $S$, $S_{1}$ is its blowing-up, and 
for $i>1$, $p_i$ is a point on $S_{i-1}$ and $S_i$ 
is its blowing-up.
We denote $p_i(K)$ the $i$-th point of an ordered cluster $K$.

Note that a point $p$ infinitely near to $q\in S$ is a
proper point on a well determined smooth surface $S(p)$, 
obtained from $S$ by blowing up all the (finitely many) 
points preceding $q$. If $C$ is a curve on $S$, let $\tilde C$
be the strict transform of $C$
in $S(p)$. We say that $p$ belongs to $C$ (as an infinitely
near point) if and only if it belongs to $\tilde C$, and
we write $\mult_p(C)=\mult_p(\tilde C)$.
Let $p, q$ be a pair of proper or infinitely near points
in $S$. The point $q$ is said to be \emph{proximate} to $p$
if and only if $p<q$ and $q$ belongs, as a proper or infinitely
near point, to the exceptional divisor
of blowing up $p$.
We shall write $p \prec q$
to mean that $q$ is proximate to $p$.

A proper point of $S$ is proximate to no other point.
An infinitely near point $p$ is always proximate to its 
immediate predecessor $\bar p$, and to at most another point $q$,
which must satisfy $q \prec \bar p$ \cite[5.3]{Cas00}. 
An infinitely near point proximate to two points
is called a \emph{satellite} point; all other points (including
proper points) are called \emph{free}.


Given a point $p$ of a cluster $K$, let $E_p$ be the
exceptional divisor of blowing up $p$, on the surface
$S_p$ obtained from $S$ by blowing up $p$ and all points
preceding it. Blowing up the remaining points gives a
morphism $\pi:S_K \longrightarrow S_p$, and we may consider
the total transform $\bar E_p$ and the strict transform
$\tilde E_p$ of $E_p$ under the composition of the blowing-ups.
We usually denote the
total transform as $E_p=\bar E_p$, if this does not lead to 
confusion. 



The set of points of a cluster $K$, equipped with the proximity
relation, has an abstract combinatorial structure, which Enriques
encoded in a convenient diagram, now called the
\emph{Enriques diagram} of the cluster 
(see \cite[IV.I]{EC15}, \cite{Cas00}).  It will
be convenient for us to give a formal definition
of Enriques diagrams along the lines of the one given by
Kleiman and Piene in \cite{KP99} (see also \cite{GSG92}):

A \emph{tree} is a finite directed graph, without
loops; it has a single initial vertex, or \emph{root}, and every other
vertex has a unique immediate predecessor.
If $p$ is the immediate predecessor of the vertex
$q$, we say that $q$ is a successor of $p$. 
An \emph{Enriques diagram} is a finite union of trees
with a binary relation between vertices, called \emph{proximity},
which satisfies:
\begin{enumerate}
\item The roots are proximate to no vertex.
\item Every vertex that is not a root is proximate to
   its immediate predecessor.
\item No vertex is proximate to more than two vertices.
\item If a vertex $q$ is proximate to two vertices
   then one of them is the immediate predecessor of $q$,
   and it is proximate to the other.
\item Given two vertices $p$, $q$ with $q$ proximate to $p$,
   there is at most one vertex proximate to both of them.
\end{enumerate}
The vertices which are proximate to two points are called
\emph{satellite}, the other vertices are called \emph{free}.
We usually denote the set of vertices of an Enriques diagram
$\Du$ with the same letter $\Du$.
We shall consider two trees with the same number of vertices
and satisfying the same proximity relations to be
the same Enriques diagram.

To show graphically the proximity relation, Enriques diagrams are
drawn according to the following rules:
\begin{enumerate}
\item If $q$ is a free successor of $p$
   then the edge going from $p$ to $q$ is smooth and curved and,
   if $p$ is not a root, it has at $p$ the same tangent
   as the edge joining $p$ to its predecessor.
\item The sequence of edges connecting a maximal
   succession of vertices proximate to the same vertex $p$
   are shaped into a line segment, orthogonal to the edge joining $p$
   to the first vertex of the sequence.
\end{enumerate}



Following \cite{KP99}, we define the following
numerical invariants of Enriques diagrams:
\begin{align*}
  \frs(\Du) &:= \text{the number of free vertices in }\Du, \\
  \rts(\Du) &:= \text{the number of roots in }\Du, \\
  \dim(\Du) &:= \rts(\Du)+\frs(\Du).
\end{align*}
$\dim (\Du)$ is called the \emph{dimension} of the Enriques
diagram, and can be thought of as the number of degrees of
freedom for the points of clusters with Enriques diagram $\Du$.

It is clear how to associate to each cluster $K$ its Enriques
diagram $\Du$: it has a vertex for each point of $K$, and its 
vertices satisfy the same proximity relations as the corresponding
points of $K$. In particular, the proper points of $K$
correspond to the roots of $\Du$. 
The converse statement is also true, namely, given an Enriques
diagram $\Du$, there exist clusters of points of $S$ whose
diagram is $\Du$ (this is implicit in \cite{KP99}).

A tree, and therefore an Enriques diagram, comes with
a natural ordering $\le$. The natural ordering of 
vertices in the Enriques diagram of a cluster $K$ 
corresponds exactly to the natural ordering of the 
points of $K$. As in the case of clusters, an admissible
ordering for an Enriques diagram $\Du$ is a total ordering
refining the natural ordering, and an ordered Enriques
diagram $\D$ is an Enriques diagram equipped with an
admissible ordering. We denote $p_i(\D)$ the $i$-th
vertex of $\D$.

  Let $K=(p_1, p_2, \dots, p_r)$ be an ordered cluster. 
The proximity matrix $P$ of $K$ is the square $r \times r$
matrix whose entry in the $i$-th row and $j$-th column is
$$
p_{i,j}=
\begin{cases}
  1 & \text{if }i=j,\\
 -1 & \text{if $p_i$ is proximate to $p_j$},\\
  0 & \text{otherwise.}
\end{cases}
$$ 

\begin{Rem}
The proximity matrix of an ordered cluster depends only
on the ordering and the proximities satisfied by its points. Therefore,
all clusters with the same Enriques diagram have the
same proximity matrix, and we speak of the
proximity matrix of an ordered Enriques diagram.
Moreover,  two ordered Enriques
diagrams are equal if and only if they have the
same proximity matrix.
\end{Rem}

It is well known that $\Pic S_K$$\cong \Pic S \oplus \bigoplus_{i=1}^r \Z E_i$,
and $E_i^2=-1$, $E_i\cdot E_j=0$ for all $i \ne j$.
Moreover,  if $D$ is an effective divisor on $S_K$, linearly 
equivalent to $\sum_{p\in K} n_pE_p$ for
some $n_p$, $p\in K$, then $D=\sum_{p\in K} n_pE_p$,
and $n_p \ge 0$ for all $p\in K$.
Note that $(E_1, E_2, \dots, E_r)$ and
$(\tilde E_1, \tilde E_2, \dots, \tilde E_r)$
are two different bases of the $\Z$--module
$\bigoplus_{i=1}^r \Z E_i$. In fact,
the invertible matrix $P$ is the matrix of
base change between the two:
\begin{equation}
\label{cbase}
 (\tilde E_1, \tilde E_2, \dots, \tilde E_r)=
   (E_1, E_2, \dots, E_r) P
\end{equation}
(see \cite[4.4]{Cas00} or \cite[1.1.29]{Alberich}).
The following two lemmas will be useful later on.

\begin{Lem}
Let $K=(p_1, p_2, \dots, p_r)$ be an ordered cluster of 
points of $S$, $P$ its proximity matrix. A divisor
$D=\sum_{i=1}^r m_iE_i$, $m_i \in \Z$, on $S_K$ is
effective if and only if
$$
P^{-1}\m \ge 0,
$$
where $\m=(m_1, m_2, \dots, m_r)$ is the column vector
of the coefficients of $D$.
\end{Lem}
\begin{proof}
See \cite[1.1.45]{Alberich}.
\end{proof}

\begin{Lem}
  \label{compartir}
Let $K=(p_1, p_2, \dots, p_r)$ be an ordered cluster of points of $S$, $P$
its proximity matrix. Two effective divisors
$D=\sum_{i=1}^r m_iE_i$, 
$D'=\sum_{i=1}^r m_i'E_i$, $m_i, m_i' \in \Z$, on $S_K$ 
have common components if and only if
$$
(\m')^T (PP^T)^{-1}\m > 0,
$$
where $\m=(m_1, m_2, \dots, m_r)$ and 
$\m'=(m_1', m_2', \dots, m_r')$ are the column vector
of the coefficients of $D$ and $D'$ respectively.
\end{Lem}
\begin{proof}
All irreducible components of $D$ and $D'$ are among
$\tilde E_1, \tilde E_2, \dots, \tilde E_r$, and because of
(\ref{cbase}), we have
\begin{align*}
D&=(E_1, E_2, \dots, E_r) \m =
(\tilde E_1, \tilde E_2, \dots, \tilde E_r) P^{-1} \m, \\
D'&=(E_1, E_2, \dots, E_r) \m' =
(\tilde E_1, \tilde E_2, \dots, \tilde E_r) P^{-1} \m'. 
\end{align*}   
As $D$ and $D'$ are effective, no coefficients in the
vectors $P^{-1} \m$, $P^{-1} \m'$ are negative, so it is
clear that $D$ and $D'$ have common components if and
only if the product $(P^{-1} \m')^T(P^{-1} \m)$ is
nonzero (and hence positive) as claimed.
\end{proof}



We shall be concerned with families
of smooth surfaces and relative divisors on them.
For our purposes, a family of surfaces is a smooth
morphism of relative dimension 2,
$$
\psi: S \longrightarrow T
$$
where $T$ is a variety defined over the algebraically
closed field $k$. $T$ is called the \emph{parameter space}
of the family, and the fibers of $\psi$, which are
smooth surfaces, are the \emph{members} of the family.
For every point $t\in T$, the fiber over $t$ will be
denoted, as customary, $S_t=S \times_T \{t\}$. As a set,
$S_t=\psi^{-1}(t)$. If $T' \rightarrow T$ is an
arbitrary morphism, then by base change we obtain a
new family of surfaces; we shall usually denote
$S_{T'}= S \times_T T'$, and the new family will be
$$
\psi_{T'}: S_{T'} \longrightarrow T' \, .
$$
A \emph{relative divisor} on the family $S \rightarrow T$
is a Cartier divisor on $S$ which meets properly every
member of the family. 

Our interest here is focused in families of surfaces
obtained by blowing up families of points or families
of clusters in surfaces or families of surfaces. Let
us now recall the construction of families of blowing-ups,
which Kleiman developed in all generality in \cite{Kle81a},
in the particular case we are dealing with.

Let $\psi:S \rightarrow T$ be a family of surfaces, 
with $T$ a smooth variety, 
and let $i:Y \hookrightarrow S$ be a closed
embedding. Consider the fiber product $S_Y=S \times_T Y$
and the diagonal morphism
$\Delta := i \times_T Id_Y: Y \longrightarrow S_Y$.
$\Delta$ is a smooth embedding over $\psi_Y$
(see \cite[17.3]{EGA4}), and its
image $\Delta(Y)$ is a closed subvariety
isomorphic to $Y$. Consider the blowing-up
$$ 
\begin{CD}
\BF(S,Y,T):=\Bl(S_Y,\Delta(Y))
     @>\pi_{\Delta(Y)}>> S_Y,
\end{CD}
$$
and the commutative diagram
$$
\diagramcompileto{familia}
    \BF(S,Y,T) \rrto^(.55){\quad p \circ \pi_{\Delta(Y)}}
               \dto_{\psi_Y \circ \pi_{\Delta(Y)}} & &S \dto \\
    Y \rrto & & T
\enddiagram
$$
We call $\pi=p \circ \pi_{\Delta(Y)}$ and
$\psi'=\psi_Y \circ \pi_{\Delta(Y)}$.
As $\Delta$ is a smooth embedding over $\psi_Y$, 
it follows that $\psi'$ is smooth of relative
dimension 2 (cf. \cite[19.4]{EGA4}), so it is
a family of smooth surfaces. We call
$$\BF(S,Y,T) \overset{\psi'}\longrightarrow Y$$
the family of blowing up $S$ at the points of $Y$.
The name is justified by the next proposition:

\begin{Pro}
\label{kleiman}
For every point $y \in Y$, and $t=\psi(y) \in T$,
consider the blowing-up $\pi_{y}: \Bl(X_t,\{y\}) \rightarrow X_t$.
Then there is a unique isomorphism
$$\Bl(X_t,\{y\}) \overset{\eta}\longrightarrow \BF(X,Y,T)_y$$
satisfying $\pi_{y}=\pi|_{\BF(X,Y,T)_y} \circ \eta$.
\end{Pro}
\begin{proof} Follows from \cite[2.4]{Kle81a}, as $\Delta(Y)$ is
obviously a local complete intersection, flat over $Y$.
\end{proof}

Let now $S$ be a smooth irreducible projective algebraic surface.
Iterating the process of blowing-up families it is possible to
define varieties parametrizing the clusters with $r$ points of $S$,
as follows. Take $X_{-1}=Spec \, k$, $X_0=S$,
$\psi_0:S \rightarrow Spec \, k$, and
define recursively $X_i$, $\psi_i$ as the blowing-up family
$$
   X_i=\BF(X_{i-1},X_{i-1},X_{i-2})
      \overset{\psi_i}\longrightarrow X_{i-1}.
$$
The morphism $\psi_i$ is
in this case projective, so its fibers are projective
smooth surfaces. For every $i$, the variety $X_i$ is
irreducible and smooth of dimension $2i+2$. The
construction of $X_i$ gives also morphisms 
$X_i \rightarrow X_{i-1}$ 
whose restrictions to the fibers of $\psi_i$ are,
by Proposition \ref{kleiman}, 
the blowing-ups of the points of the fibers of $\psi_{i-1}$;
we denote these morphisms $\pi_i:X_i \rightarrow X_{i-1}$.
To simplify notations, let us put
$\pi_{r,i}=\pi_{i+1} \circ \pi_{i+2} \circ \cdots \circ \pi_r$,
$\psi_{r,i}=\psi_i \circ \psi_{i+1} \circ \cdots \circ \psi_{r}$.
For any point $x \in X_i$, we call
$$S(x)=(X_i)_{\psi_i(x)}$$
the blown up surface containing $x$.
Recall that for any cluster $K$, 
$\pi_K:S_K \rightarrow S$ is the composition
of the blowing-ups of the points in $K$.

The following proposition makes the set of all
ordered clusters with $r$ points into an
algebraic variety.

\begin{Pro}
\label{clustersfam}
For every $r \geq 1$ there is a bijection
$$ X_{r-1} \overset{K}\longrightarrow
   \left\{ \text{\emph{ordered clusters of }}
            r \text{ \emph{points}} \right\}
$$
and, for every $x\in X_{r-1}$, a unique isomorphism
$\eta_x:S_{K(x)}\rightarrow (X_r)_x$ such that
\begin{enumerate}
\item \label{p1} $\pi_{K(x)}=\pi_{r,0}|_{(X_r)_x} \circ \eta_x$.
\item \label{p2} If $K(x)=(p_1, p_2, \dots, p_r)$ then
$K(\psi_{r-1,i}(x))=(p_1, p_2, \dots, p_i)$, i.e.
$\psi_{r-1,i}$ maps the point of $X_{r-1}$ corresponding to
a cluster to the point of $X_{i-1}$ corresponding to the
cluster of its first $i$ points.
\item \label{p3} If $K(x)=(p_1, p_2, \dots, p_r)$ then, for every 
$p \in S_{K(x)}$, $K(\eta_x(p))=(p_1, p_2, \dots, p_r, p)$.
\end{enumerate}
\end{Pro}

It is interesting to note that the existence
of a bijection $K$ follows from \cite[I.2]{Har88},
thanks to the obvious bijection
\begin{align*}
 \left\{ \text{\emph{ordered blowing-ups at }}
            r \text{ \emph{points}} \right\}
                    &  \longrightarrow
 \left\{ \text{\emph{ordered clusters of }}
            r \text{ \emph{points}} \right\} \\
 S_K &\longmapsto K
\end{align*}
However, the proof we are now going to give has the 
advantadge of being more explicit, and we also
obtain properties \ref{p1}, \ref{p2} and \ref{p3}
of the isomorphisms $\eta_x$, which will be
very useful for our study of $X_r$.
Notice also that the ordering of points in clusters
is essential in Proposition \ref{clustersfam}.
Using unordered clusters would give rise to 
a non-injective $K$.

\begin{proof}
First of all we define the map $K$. Given
$x=x_r\in X_{r-1}$, consider the points
$x_i=\psi_i(x_{i+1}) \in X_{i-1}$ for $i\in \{r-1, r-2, \ldots, 1\}$.
Proposition \ref{kleiman} applied to
$X_i=\BF(X_{i-1},X_{i-1},X_{i-2})$ shows that
there are unique isomorphisms
$$\Bl(S(x_i),\{x_i\})
   \overset{\eta_{x_i}}\longrightarrow (X_i)_{x_i}$$
such that $\pi_i|_{(X_i)_{x_i}} \circ \eta_{x_i}$
is the blowing-up of $x_i$,
and obviously $x_{i+1}\in S(x_{i+1})=(X_i)_{x_i}$ for all $i$.
So $K(x):=(x_1, \eta_{x_1}^{-1} (x_2), \ldots, \eta_{x_{r-1}}^{-1} (x_r))$
is a cluster and $\eta_{x_r}$ is the claimed isomorphism
(i.e.,  properties \ref{p1}, \ref{p2} and \ref{p3} are satisfied).

It remains to be seen that $K$ is bijective.
Let $K_0=(p_1, p_2, \ldots, p_r)$ be a cluster
with $r$ points. We have $p_1 \in S=X_0$,
$p_2 \in \Bl(X_0,\{p_1\}) \cong (X_1)_{p_1} \subset X_1$,
and iterating the process $r$ times,
$$
\begin{CD}
p_r \in \Bl(S(p_{r-1}),\{p_{r-1}\})
   @>{\eta_{p_{r-1}}}>{\cong}> 
   (X_{r-1})_{\eta_{p_{r-2}}(p_{r-1})} \subset X_{r-1}.
\end{CD}
$$
So we have $\eta_{p_{r-1}}(p_r)\in X_{r-1}$, and it is clear
that $K(\eta_{p_{r-1}}(p_r))=K_0$, so we have seen
that the map $K$ is onto.

To see that $K$ is injective, we
will use induction on $r$. For $r=1$, the claim
is obvious. For $r>1$, suppose $K(x)=K(x')$, with
$x, x' \in X_{r-1}$.
Then in particular $K(\psi_{r-1}(x))=K(\psi_{r-1}(x'))=K_1$,
by the definition of $K$, and we have isomorphisms
\begin{gather*}
S_{K_1} \overset{\eta}\longrightarrow S(x) \\
S_{K_1} \overset{\eta'}\longrightarrow S(x'),
\end{gather*}
with
$\pi_{K_1}=\pi_{r-1,0}|_{S(x)} \circ \eta
          =\pi_{r-1,0}|_{S(x')} \circ \eta'$.
Both clusters being equal, their last point must
be the same, that is,
$\eta^{-1}(x)={\eta'}^{-1}(x')$.
On the other hand, the induction hypothesis
says that $\psi_{r-1}(x)=\psi_{r-1}(x')=x_{r-1}$, so
$S(x)=S(x')=(X_{r-1})_{x_{r-1}}$, and
$\eta=\eta'$. As $\eta$ is an isomorphism,
$\eta^{-1}(x)=\eta^{-1}(x')$ implies $x=x'$.
\end{proof}




\section{Fixing the Enriques diagram}
\label{fix}

From now on we identify the set of ordered clusters with $r$ 
points to the variety $X_{r-1}$, and denote the
points of $X_{r-1}$ by capital letters such as $K$.
Thus, given an ordered Enriques diagram $\D$ of
$r$ vertices, the clusters with Enriques diagram $\D$ 
describe a subset $\Cl(\D) \subset X_{r-1}$.
Our next goal is to see that $\Cl(\D)$ is a locally
closed smooth subvariety of $X_{r-1}$.

Throughout this section, we
work on the family of surfaces
$$
\psi_r:X_r \longrightarrow X_{r-1}
$$
and certain subfamilies of this one; to lighten a little
bit the notations, for every subvariety $Y\subset X_{r-1}$
we shall write $S_Y=(X_r)_Y$ and let 
$\psi_r|_{S_Y}:S_Y \longrightarrow Y$
be the family of surfaces corresponding to the set of
clusters $Y$. For example,
once we prove that $\Cl(\D)$ is a locally closed smooth
subvariety of $X_{r-1}$, $S_{\Cl(\D)}\rightarrow \Cl(\D)$
will be the family of all surfaces obtained by blowing
up clusters whose Enriques diagram is $\D$.



The proximity relations which are encoded in
the Enriques diagram of a cluster express the belonging
of some points of the cluster to exceptional divisors of
previously blown-up points. To deal with $\Cl(\D)$,
it will be necessary to have some knowledge of the
\emph{families} formed by these divisors in the family of
the blown-up surfaces.
Call $F_i$ the pullback by $\pi_{r,i}$ in $X_r$ of the
exceptional divisor of the $i$-th blowing up,
$$ \begin{CD}
X_i @>{\pi_{\Delta(X_{i-1})}}>> X_{i-1} \times_{X_{i-2}} X_{i-1}.
\end{CD} $$
It follows from Proposition \ref{kleiman} that,
for any $K=(p_1, p_2, \dots, p_r) \in X_{r-1}$, the
restriction of $F_i$ to $S_K \subset X_r$ is the pullback $E_i$ of 
the exceptional divisor of blowing up $p_i$.
Therefore, every $F_i$ is an effective relative divisor over $X_{r-1}$, and
moreover it is irreducible (if the restriction of an effective 
relative divisor to some member of the family is irreducible,
and the base of the family is irreducible, then the relative
divisor is irreducible also). 

Given an ordered Enriques diagram $\D$, let $P$ be its proximity matrix.
Recall that $\tilde E_i$ denotes the \emph{strict}
transform in $S_K$ of the exceptional divisor of blowing up $p_i$,
and a cluster $K$ has Enriques diagram $\D$
if and only if
$$ (\tilde E_1, \tilde E_2, \dots, \tilde E_r)=
   (E_1, E_2, \dots, E_r) P$$
in $S_K$.
For any cluster $K$ with $r$ points, we define the
\emph{virtual exceptional divisors} in $S_K$ relative
to the Enriques diagram $\D$ to be
$$ (E_1^{\D}, E_2^{\D}, \dots, E_r^{\D})=
   (E_1, E_2, \dots, E_r) P \, .$$
The first step toward the construction of $\Cl(\D)$ is to consider the 
set $\Eff(\D)$ of all clusters $K$ for which the virtual exceptional
divisors $E_1^{\D}, E_2^{\D}, \dots, E_r^{\D}$ are effective in $S_K$.
We shall see that $\Eff(\D)$ is a closed subvariety
of $X_{r-1}$ and $\Cl(\D)$ is open in $\Eff(\D)$. This
will prove that $\Cl(\D)$ is locally closed.

\begin{Rem}
The definition of $(E_1^{\D}, E_2^{\D}, \dots, E_r^{\D})$
immediately gives that the Enriques diagram of $K$ is $\D$
if and only if $\tilde E_i = E_i^{\D} \ \forall i$.
In particular, $\Cl(\D) \subset \Eff(\D)$, because
the divisors $\tilde E_i$ are always effective.
\end{Rem}

Define divisors $F_i^{\D}$ on $X_r$ as
$$ (F_1^{\D}, F_2^{\D}, \dots, F_r^{\D})=
   (F_1, F_2, \dots, F_r) P . $$
For every subvariety $Y\subset X_{r-1}$, we shall
write $F_{i,Y}^{\D}$ the restriction of $F_i^{\D}$
to $S_Y$. Note that the divisors 
$F_i^{\D}$ are relative divisors over $X_{r-1}$,
therefore $F_{i,Y}^{\D}$ 
are also relative divisors over $Y$.

\begin{Pro}
\label{efclo}
The subset $\Eff(\D)$ is closed in $X_{r-1}$  
\end{Pro}

\begin{proof}
Proposition \ref{kleiman} tells us that the restriction
of $F_i^{\D}$ to any surface $S_K \subset X_r$ is exactly
$E_i^{\D}$. In other words, $F_{i,K}^{\D}=E_i^{\D}$, or 
$(\O_{X_r}(F_i^{\D}))_K=\O_{S_K}(E_i^{\D})$. So
$\Eff(\D)$ can be described as
$$
\Eff(\D)=\{K \in X_{r-1} \,|\, h^0(S_K,(\O_{X_r}(F_i^{\D}))_K) > 0 \}\, 
$$
($h^0(S_K,\O_{S_K}(E_i^{\D})) > 0$ means that there
is an effective divisor $D$ in $S_K$ linearly equivalent to
$E_i^{\D}$; 
this implies
$D=E_i^{\D}$ and therefore $E_i^{\D}$ itself is effective).
On the other hand, $X_r \rightarrow X_{r-1}$ is smooth, so
the invertible sheaf $\O_{X_r}(F_i^{\D})$ is flat
over $X_{r-1}$, and by the semicontinuity theorem
\cite[III,12.8]{HAG} the claim follows. 
\end{proof}

We want to see next that $\Cl(\D)$ is open in $\Eff(\D)$; for that
we shall give some attention to the complementary set
$\Eff(\D) \setminus \Cl(\D)$, that is, to the clusters of $\Eff(\D)$
wich have Enriques diagrams other than $\D$.

\begin{Lem}
\label{noiguals}
\label{incidencia}
Let $\D$, $\D'$ be two Enriques diagrams with $r$
vertices whose proximity matrices are $P$, $P'$, and
let $K$ be an ordered cluster whose
Enriques diagram is $\D'$. The following are equivalent:
\begin{enumerate}
\item \label{unc} $K \in \Eff(\D)$.
\item \label{totc} $\Eff(\D') \subset \Eff(\D)$.
\item \label{mat} All entries in the matrix ${P'}^{-1}P$ are
non negative.
\end{enumerate}
Moreover, $\Eff(\D') = \Eff(\D)$ if and only if $\D'=\D$.
\end{Lem}
\begin{proof}
By definition of the proximity matrices we have
\begin{gather*}
 (\tilde E_1, \tilde E_2, \dots, \tilde E_r)=
   (E_1, E_2, \dots, E_r) P' \\
 (E_1^{\D}, E_2^{\D}, \dots, E_r^{\D})=
   (E_1, E_2, \dots, E_r) P \ 
\end{gather*}   
on $S_{K}$, 
therefore $(E_1^{\D}, E_2^{\D}, \dots, E_r^{\D})=
   (\tilde E_1, \tilde E_2, \dots, \tilde E_r) {P'}^{-1}P$,
so the entries in the matrix ${P'}^{-1}P$ are the coefficients
of the expression of each $E_i^{\D}$ as a linear combination of 
irreducible divisors. All divisors $E_i^{\D}$ are effective if and
only if these coefficients are nonnegative, 
and the equivalence of \ref{unc} and \ref{mat} follows. 

To see the equivalence with \ref{totc}, it will be enough
to prove that \ref{mat} implies \ref{totc}, because
clearly \ref{totc} implies \ref{unc}. So assume that 
all entries in the matrix ${P'}^{-1}P$ are
non negative, and let $K'$ be a cluster in $\Eff(\D')$;
we have to see that $K' \in \Eff(\D)$. In $S_{K'}$, it holds
$(E_1^{\D}, E_2^{\D}, \dots, E_r^{\D})=
   (E_1^{\D'}, E_2^{\D'}, \dots, E_r^{\D'}) {P'}^{-1}P.$
Therefore, all divisors $E_i^{\D}$ are linear combinations
with positive coefficients
of the divisors $E_i^{\D'}$, which are effective 
(but may be reducible!) because
$K \in \Eff(\D')$, so they are effective and $K' \in \Eff(\D)$.

It remains to be seen that $\D' \neq \D$ implies
$\Eff(\D') \ne \Eff(\D)$.
Let $i$ be the least index of a vertex in which the
two diagrams differ, i.e., the proximities satisfied
by the vertices $p_j(\D)$ and $p_j(\D')$ are the
same for $j<i$, but not for $j=i$.
This means that there is $j<i$ with
either $p_j(\D) \prec p_i(\D)$, $p_j(\D') \nprec p_i(\D')$ or
$p_j(\D) \nprec p_i(\D)$, $p_j(\D') \prec p_i(\D')$; both
cases are symmetric, so assume
$p_j(\D) \nprec p_i(\D)$, $p_j(\D') \prec p_i(\D')$.
For every cluster $K$
whose Enriques diagram is $\D$, $E_j^{\D'}$ is not
effective in $S_K$, therefore $K\not\in \Eff(\D')$
and $\Eff(\D') \neq \Eff(\D)$.
\end{proof}

\begin{Pro}
\label{loclo}
$\Cl(\D)$ is a nonempty open subset of $\Eff(\D)$.  
\end{Pro}
\begin{proof}
By definition, $\Cl(\D)$ is nonempty.
Lemma \ref{noiguals} implies that
$$
\Cl(\D)=\Eff(\D) \setminus \bigcup_{\Eff(\D')\varsubsetneq \Eff(\D)} \Eff(\D') \, .
$$
Remark that the union is finite, because there are just a finite
number of Enriques diagrams with $r$ vertices.
The variety $\Eff(\D)$ being closed for all $\D$, we can conclude
that $\Cl(\D)$ is a nonempty open subset of $\Eff(\D)$. 
\end{proof}

Proposition \ref{loclo} shows that the variety
$X_{r-1}$ decomposes as a disjoint union of locally closed subsets 
corresponding to all possible Enriques diagrams with $r$ points:
$$X_{r-1} = \bigcup_{\sharp \D=r} \, \Cl(\D) \, .$$
It is interesting to note that, as it will be proved in
section \ref{nostrata}, this decomposition is \emph{not}
a stratification, i.e., it is not true that
$\Cl(\D')\cap \overline {Cl(\D)}\ne \emptyset$
implies $\Cl(\D') \subset \overline {\Cl(\D)}$.

We have shown
that $\Eff(\D)$ is closed in $X_r$ using
the semicontinuity theorem.
One can also give a more explicit inductive construction.
If $\D$ is the diagram consisting of a single point, then
$\Eff(\D)=X_0=S$. If $\D$ has more than one vertex, let 
$\smash{\breve \D}$ be the Enriques diagram
obtained from $\D$ by dropping the last vertex, and
for any $K$, let $\smash{\breve K}$ be
the cluster obtained from $K$ by dropping the last point.
Consider the family of surfaces
$$S_{\Eff(\breve \D)} \rightarrow \Eff(\breve \D) \, .$$
Its fibers are exactly the 
$S_{\breve K}$ which have $E_i^{\breve \D}$ effective, $i=1, \dots, r-1$.

Recall from Proposition \ref{clustersfam} that the point of 
$X_{r-1}$ to which we identify the
cluster $K$ is exactly
$\smash{p_r(K) \in S_{\breve K} \subset X_{r-1}}$; we also have
$\psi_{r-1}(K)=\breve K$.
As $K\in \Eff(\D)$ trivially implies $\breve K \in \Eff(\breve \D)$,
it follows that
$$
\Eff(\D) \subset S_{\Eff(\breve \D)} \, .
$$
Let
$$
\pi: S_K \longrightarrow S_{\breve K}
$$
be the blowing-up of the last point, $p_r(K)$. Then,
by its definition,
$$
E_i^{\D}=\begin{cases}
 E_r & \text { if } i=r \\
 \pi^{*}(E_i^{\breve \D}) & \text{ if } i<r \text{ and }
         p_i(\D) \not\prec p_r(\D) \\
 \pi^{*}(E_i^{\breve \D}) -E_r & \text{ if } i<r \text{ and }
         p_i(\D) \prec p_r(\D) \, .
\end{cases}  
$$
Assuming $\breve K \in \Eff(\breve\D)$, these divisors will
be effective if and only if
$$
 p_r(K) \in E_i^{\breve \D} \ \ \forall i \text{ s.t. } p_i(\D) \prec p_r(\D) \, .
$$
This allows us to describe $\Eff(\D)$ precisely:
we have in fact proved

\begin{Pro}
\label{construccioXD}
Let $F_i(\breve \D)$ be the restriction of $F_i^{\breve \D}$ to
$S_{\Eff(\breve \D)}$; then
$F_i(\breve \D)$ is effective for $i=1, 2, \dots, r-1$ and
$$\Eff(\D)=
\begin{cases}
 S_{\Eff(\breve \D)} & \negthickspace \text{if } p_r(\D) \text{ is a root,} \\
 F_i(\breve \D) & \negthickspace \text{if } p_r(\D) \text{ is free, } 
          p_r(\D) \succ p_i(\D) \\
 F_i(\breve \D) \cap F_j(\breve \D) & \negthickspace \text{if } p_r(\D) 
     \text{ is satellite, } p_r(\D) \succ p_i(\D), p_j(\D), \ i \ne j . 
\end{cases}  
$$
\end{Pro}

\begin{Pro}
\label{localtancat}
$\Cl(\D)$ is irreducible and smooth of dimension $\dim \D$.
\end{Pro}

\begin{proof}
By induction on the number of points $r$
of the cluster. For $r=1$ there is only one possible
Enriques diagram, and $\Cl(\D)=X_0=S$ is irreducible
and smooth of dimension 2.
For $r>1$ consider the Enriques diagram $\breve \D$
obtained from $\D$ by dropping the last point. It
is clear that $\psi_{r-1}(\Cl(\D)) =\Cl(\breve \D)$.
The morphism
$$ \begin{CD}
\Cl(\D) @>\psi_{r-1}>> \Cl(\breve \D)
\end{CD} $$
is smooth because it comes from 
$\psi_{r-1}:X_{r-1} \rightarrow X_{r-2}$
by base change, 
and has irreducible fibers (the fibers are
single points if the last point of $\D$ is satellite,
open sets of $\P^1$ if it is free, and open sets of
$S$ if it is a root). The irreducibility and smoothness
of $\Cl(\D)$ then follow from those of $\Cl(\breve \D)$, which
we obtain from the induction hypothesis; adding up
the dimension of the fibers to
$\dim \Cl(\breve \D) = \dim \breve \D$
we obtain also the claimed dimension of $\Cl(\D)$.
\end{proof}

At this point we introduce the notion of specialization of Enriques
diagrams. We are interested in the incidence relations between
the varieties corresponding to different Enriques diagrams.
To be precise, we would like to know the closure
$\overline {\Cl(\D)}$ of $\Cl(\D)$ in $X_{r-1}$, and in particular
to know for which couples of diagrams $\D, \D'$
there is an inclusion
$\Cl(\D') \subset \overline{\Cl(\D)}$.
We will henceforth say that the ordered Enriques diagram $\D$
\emph{specializes to} $\D'$ or that $\D'$ 
\emph{is a specialization of} $\D$,
$\D \rightsquigarrow \D'$, whenever 
$\Cl(\D') \subset \overline{\Cl(\D)}$.
After \ref{noiguals}, it is easy to see that
specialization is an order relation on the set of
Enriques diagrams with $r$ vertices, for each $r>0$. It is not
a total order, as shown by example \ref{noespe} below.

\begin{Lem}
\label{necesp}
  Let $\D$, $\D'$ be two ordered Enriques diagrams
whose proximity matrices are $P$ and $P'$ respectively,
and assume $\D \rightsquigarrow \D'$. Then ${P'}^{-1}P$
has no negative entries.
\end{Lem}
\begin{proof}
  $\Eff(\D)$ being closed (\ref{efclo})
and $\Cl(\D)\subset \Eff(\D)$, it follows that
$\overline{\Cl(\D)}\subset \Eff(\D)$, so the
hypothesis $\D \rightsquigarrow \D'$ says
$\Cl(\D') \subset \Eff(\D)$. Now
\ref{noiguals} gives the claim.
\end{proof}

\begin{Exa}
\label{noespe}
Consider the two ordered Enriques diagrams
of figure \ref{noespfig}, and let their
proximity matrices be $P$ and $P'$. It is
immediate to see that both ${P'}^{-1}P$
and $P^{-1}P'$ have negative entries,
so lemma \ref{necesp} proves that neither of
them is a specialization of the other.
 \begin{figure}
  \begin{center}
    \mbox{\includegraphics{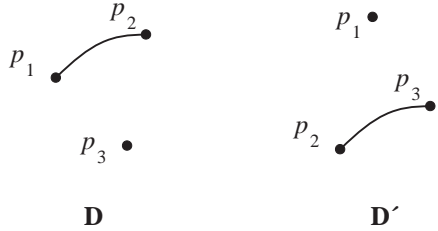}}
    \caption{Specialization is not a total order:
$\D \not\rightsquigarrow \D'$, $\D' \not\rightsquigarrow \D$.}
    \label{noespfig}
  \end{center}
\end{figure}
\end{Exa}

\begin{Lem}
  \label{efcomp}
Every irreducible component of the variety $\Eff(\D)$
has the form $\overline{\Cl(\D')}$ for some Enriques
diagram $\D'$ such that $\Eff(\D') \subset \Eff(\D)$. Moreover,
one of the components of $\Eff(\D)$ is $\overline{\Cl(\D)}$.
\end{Lem}
\begin{proof}
Because of \ref{noiguals}, and the closedness of
$\Eff(\D)$, we have
$$ \Eff(\D) = \bigcup_{\Eff(\D')\subset \Eff(\D)} {\Cl(\D')}
= \bigcup_{\Eff(\D')\subset \Eff(\D)} \overline{\Cl(\D')} \ .$$
Since the varieties $\Cl(\D')$ are irreducible, 
and the union is finite, every irreducible
component of $\Eff(\D)$ has the form $\overline {\Cl(\D')}$ for
some $\D'$ with $\Eff(\D')\subset \Eff(\D)$.
Moreover, one of these components must be $\overline {\Cl(\D)}$;
otherwise we would have
$$\Cl(\D) \subset \overline {\Cl(\D')} \subset \Eff(\D')$$
which by \ref{noiguals} implies $\Eff(\D)\subset \Eff(\D')$.
As we also have $\Eff(\D')\subset \Eff(\D)$ it follows
that $\Eff(\D) = \Eff(\D')$ 
against lemma \ref{noiguals}.

\end{proof}

$\Eff(\D)$ is not always irreducible nor even equidimensional:
it may happen that some $\D'$ with $\Eff(\D')\subset \Eff(\D)$
has dimension bigger than $\dim \D$. 
However, in many relevant cases
we shall show an equality $\overline {\Cl(\D)} = \Eff(\D)$,
which by \ref{efcomp} is equivalent to the irreducibility 
of $\Eff(\D)$.

\begin{Pro}
\label{prime}
Let $\D$ be an ordered Enriques diagram. The following
conditions are equivalent:
\begin{enumerate}
\item \label{s2} For every Enriques diagram $\D'\! \ne \! \D$ such that 
 $\Eff(\D') \subset \Eff(\D)$, $\dim \D' \! < \dim \D$.
\item \label{s0} For every Enriques diagram $\D'\! \ne \! \D$ such
  that ${P'}^{-1}P$ has no negative entries, $\dim \D' \! < \dim \D$.
\item \label{s1} For every cluster $K \in \Eff(\D)$ with Enriques 
 diagram $\D' \ne \D$, $\dim \D' < \dim \D$.
\item \label{s3} $\overline {\Cl(\D)}=\Eff(\D)$.
\item \label{s4} $\Eff(\D)$ is irreducible.
\end{enumerate}
\end{Pro}
An ordered Enriques diagram satisfying one (and therefore all)
of the preceding conditions will be called a \emph{prime} 
Enriques diagram. Remark that if $\D$ is prime then 
$\breve \D$ is prime also, by condition \ref{s4}, for example.

\begin{proof}
The equivalence of \ref{s2}, \ref{s0} and \ref{s1} is clear after
lemmas \ref{noiguals} and \ref{necesp}. Since 
$\overline {\Cl(\D)}$ is one component of $\Eff(\D)$
(\ref{efcomp}), \ref{s3} and \ref{s4} are equivalent.
To see that \ref{s3} implies \ref{s2}
is easy: if \ref{s2} is false
then there exists $\D' \ne \D$ with
$\Cl(\D') \subset \Eff(\D') \subset \Eff(\D)$ and
$$ \dim \Cl(\D') = \dim \D' \ge \dim \D = \dim \Cl(\D) ;$$
therefore $\Cl(\D') \not\subset \overline {\Cl(\D)}$ 
(recall that $\Cl(\D)$ and $\Cl(\D')$ are disjoint) and
$\Eff(\D) \not\subset \overline {\Cl(\D)}$, against \ref{s3}.

To prove that \ref{s2} implies \ref{s3},
we proceed by induction on the number 
$r$ of points of $\D$, noting that if $r>1$ and \ref{s2} is
true for $\D$, then it is also true for $\breve \D$.   
For $r=1$ we already noticed that there is only one Enriques diagram.
For $r>1$ we first recall that the divisors 
$\smash{F_{i,\Eff(\breve \D)}^{\breve \D}}$ 
are relative divisors over $\Eff(\breve \D)$.
This implies that every one of their components dominates $\Eff(\breve \D)$, 
because $\Eff(\breve \D)$ is irreducible, due to the induction hypothesis,
This implies that the divisors
$\smash{F_{i,\Eff(\breve \D)}^{\breve \D}}$ 
have no component in common. Indeed,
a common component of two such divisors would imply a common 
component of two divisors among the $E_i^{\breve \D}$ in all
$S_K, K \in \Eff(\breve \D)$; on the other 
hand for every cluster $K \in \Cl(\D) \subset \Eff(\D)$ the divisors 
$\smash{E_i^{\breve \D}}$ are irreducible, so they share no component.
This observation together with
\eqref{construccioXD} show that $\Eff(\D)$ is equidimensional.
We already know that every component of $\Eff(\D)$ is of the
form $\overline{\Cl(\D')}$ with $\Eff(\D') \subset \Eff(\D)$; 
as they must have dimension $\dim \D$, hypothesis 
\ref{s2} implies that there is only one component, 
$\overline {\Cl(\D)}=\Eff(\D)$. 
\end{proof}

\begin{Cor}
\label{espmatrius}
Let $\D, \D'$ be Enriques diagrams whose proximity matrices 
are $P, P'$, respectively. Suppose that $\D$ is prime.
Then $\D \rightsquigarrow \D'$ if and only if all entries
of the matrix of the matrix ${P'}^{-1}P$ are non negative.
\end{Cor}
\begin{proof}
As $\D$ is prime, $\Eff(\D)=\overline{\Cl(\D)}$ and 
\ref{incidencia} gives the result.
\end{proof}

\begin{Rem}
\label{roots}
Suppose that $K\in \Eff(\D)$ and $p_i(\D)\prec p_j(\D)$.
Then $p_j(K)$ belongs to the virtual exceptional divisor
of blowing up $p_i(K)$ (see \eqref{construccioXD}) and 
therefore it is not a root. So, given
two Enriques diagrams $\D$ and $\D'$ such that
$\Eff(\D') \subset \Eff(\D)$, if $p_i(\D')$ is a root
then $p_i(\D)$ must be a root too. This obviously
implies that
$$
roots(\D) \ge roots(\D').
$$
As the dimension of an Enriques diagram is 
$$
\dim \D = roots(\D)+ frees(\D),
$$
we conclude that a non-prime Enriques diagram $\D$ must have
at least one satellite vertex $p_i(\D)$ such that, for some
$\D'$, $\Eff(\D') \subset \Eff(\D)$ and $p_i(\D')$ is free.
Otherwise we would have that $p_i(\D)$ is a root, free
or satellite if and only if $p_i(\D')$ is a root, free
or satellite respectively, and it is not difficult to
see that this implies $\D=\D'$, using \eqref{construccioXD}
inductively.
In particular, all Enriques diagrams without satellites
are prime.
\end{Rem}

\subsection{Prime Enriques diagrams}

It would be interesting to have some criterion to decide whether
an ordered Enriques diagram is prime or not. Remark \ref{roots} 
shows that the fact of $\D$ being non-prime is related to the existence of
satellites $p_i(\D)$ and diagrams $\D'$ such that
$\Eff(\D') \subset \Eff(\D)$ and $p_i(\D')$ is free.
We shall therefore study under which conditions (on $\D$) 
this is possible.
The structure of the virtual exceptional divisors $E_i^{\D}$
will give us some information on that;
let us first show two different examples of non-prime Enriques diagrams
to give an idea of the kind of phenomena involved.

\begin{Exa} 
\label{ex1}
Consider the Enriques diagrams with six points of
figure \ref{doscusps}.
\begin{figure}
  \begin{center}
    \mbox{\includegraphics{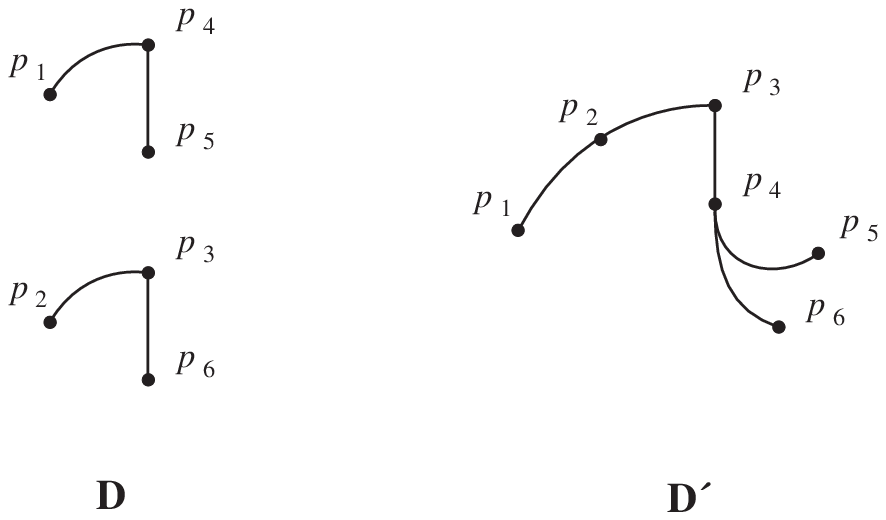}}
    \caption{$\overline{\Cl(\D')}$ is a component of $\Eff(\D)$.}
    \label{doscusps}
  \end{center}
\end{figure}
It is immediate that $\dim \D= \dim \D' = 6$,
and an easy calculation with the proximity matrices shows that
$\Eff(\D')\subset \Eff(\D)$, so $\D$ is not prime. The fifth and sixth
points are responsible for the relatively big dimension of $\D'$
compared to $\D$, because they are satellites for $\D$ and free for $\D'$.
To see why is this possible, we should look at the surface obtained
by blowing up the first four points of a cluster $K$ whose Enriques
diagram is $\D'$: $\tilde E_4$ is a component of the four divisors
$E_i^{\D}$, $i=1, \dots, 4$, so it is enough for $p_5$ 
and $p_6$ to lie in $\tilde E_4$
for all six $E_i^{\D}$ to be effective in $S_K$.
\end{Exa}

The following example shows that $\Eff(\D)$ needs
not to be of the same dimension as $\Cl(\D)$.

\begin{Exa}
\label{ex2} 
Consider the Enriques diagrams with nine points of
figure \ref{trescusps}.
\begin{figure}
  \begin{center}
    \mbox{\includegraphics{nosimpl2.eps}}
    \caption{$\overline{\Cl(\D')}$ is a component of $\Eff(\D)$.}
    \label{trescusps}
  \end{center}
\end{figure}
Now the dimensions are
$\dim \D= 8< \dim \D' = 9$,
and again an easy calculation with the proximity matrices shows that
$\Eff(\D')\subset \Eff(\D)$, so $\D$ is not prime. Now the satellites
of $\D$ which are free in $\D'$ are the fifth, seventh and nineth
points.
To see why is this possible, we should look at the surface obtained
by blowing up the first three points of a cluster $K$ whose Enriques
diagram is $\D'$: $\tilde E_3$ is a component of multiplicity 2 
of $E_1^{\D}$, so if $p_4$, $p_6$ and $p_8$ lie in 
$\tilde E_3$ then $\tilde E_4$, $\tilde E_6$ and $\tilde E_8$
will still be components of $E_1^{\D}$, so that it is enough 
for $p_5$, $p_7$ and $p_9$ to lie in $\tilde E_4$, $\tilde E_6$ and 
$\tilde E_8$ respectively in order to have $E_i^{\D}$ effective in
$S_K$, $i=1, 2, \dots, 9$.
It is not difficult to extend this example to obtain couples of
Enriques diagrams with $\Eff(\D')\varsubsetneq \Eff(\D)$ and
$\dim \D' - \dim \D$ arbitrarily big, always exploiting the double
divisor which appears in $S_3$.
\end{Exa}

Let $\D$ be an ordered Enriques diagram $\D$ with $r$ vertices.
We denote by $\D^i$ the Enriques
diagram of its first $i$ vertices (i.e. an ordered Enriques
diagram of $i$ vertices which satisfy the same proximity 
relations as the first $i$ vertices of $\D$).
Similarly, for a cluster $K$ we shall denote by $K^i$ the
cluster of the first $i$ points of $K$.

From the two examples above one can see that 
the existence of different virtual exceptional
divisors sharing components can cause the non-primality
of an Enriques diagram. More explicitly, if $p_i(\D)$ is satellite, proximate
to $p_j(\D)$ and $p_k(\D)$, and there is a cluster
$K^{i-1} \in \Eff(\D^{i-1})$ such that $E_j^{\D^{i-1}}$ and 
$E_k^{\D^{i-1}}$ share a component, then 
it is enough for $p_i \in K$ to lie in this common component
for $E_i^{\D}$ to be effective in $S_K$. So, in this conditions
\emph{there is} an Enriques diagram $\D'$  with
$\Eff(\D') \subset \Eff(\D)$ and $p_i(\D')$ free. We also 
see from the second example that virtual exceptional divisors with
multiple components easily give rise to components common to several
other virtual divisors, so they can cause non-primality of Enriques 
diagrams too.

\begin{Pro}
\label{capsat}
Let $\D$ be an ordered Enriques diagram $\D$ with $r$ vertices.
Suppose that for some $i$, $\D^i$ is prime, and no vertex
$p_j(\D)$, $j>i$, is a satellite. Then $\D$ is prime.
\end{Pro}
\begin{proof}
By induction on $r-i$. For $r=i$ there is nothing to
prove. For $r>i$ we can assume (induction hypothesis) that
$\breve \D$ is prime. Let $\D'\ne \D$ be such that 
$\Eff(\D')\subset \Eff(\D)$; we have to prove that 
$\dim \D'< \dim \D$. The hypothesis on $\D$ says
that $p_r(\D)$ is no satellite, so remark \ref{roots} says
that $\dim \D- \dim \D' \ge \dim \breve \D - \dim \breve \D'$.
If $\breve \D \ne \breve \D'$ then 
$\dim \breve \D - \dim \breve \D'>0$, because $\breve \D$
is prime, and we are done. If $\breve \D= \breve \D'$,
then the hypothesis $\D \ne \D'$ implies that $p_r(\D')$
satisfies more proximity relations than $p_r(\D)$ and 
therefore we also get $\dim \D - \dim \D'>0$.
\end{proof}

\begin{Pro}
Let $\D$ be an Enriques diagram and $K \in \Eff(\D)$ a cluster
whose Enriques diagram is $\D' \ne \D$. Consider the (effective)
divisors $E_i^{\D} \subset S_K, \ i=1, 2, \dots, r$.
\begin{enumerate}
  \item $E_i^{\D}$ is connected for all $i$.
  \item For $i \ne j$, the intersection $E_i^{\D} \cap E_j^{\D}$
is either empty, or a single point or a connected divisor.
\end{enumerate}
\end{Pro}
\begin{proof}
\begin{enumerate}
  \item By induction on the number of points $r$ of the diagram.
      For $r=1$ there is nothing to prove. For $r>1$, let $\breve \D$,
      $\breve \D'$ and $\breve K$ be the Enriques diagrams and the
      cluster obtained by dropping the last point. The induction
      hypothesis tells us that the divisors
      $E_i^{\breve \D}$ are connected. Whenever $p_i(\D) \nprec p_r(\D)$,
      the total transform of $E_i^{\breve \D}$ in $S_K$ coincides with
      $E_i^{\D}$ thus proving its connectedness. If 
      $p_i(\D) \prec p_r(\D)$, then either $p_r(K)$ is a singular point
      of $E_i^{\breve \D}$ and $E_r$ is a component of $E_i^{\D}$
      thus ensuring its connectedness or a nonsingular point of
      $E_i^{\breve \D}$ and the strict transform is connected.
  
  \item We only have to discard the possibility
      that two such divisors meet in more than one point,
      or in a divisor with more than one connected component. 
      Since we know that each divisor is connected, both possibilities
      would imply that the total exceptional divisor in $S_K$
      of a point of $S$ is not simply connected, and it is known
      that this is not possible (see \cite[4.4]{Cas00}, for example).
\end{enumerate}
\end{proof}

\begin{Pro}
\label{unsat}
If $\D$ is an ordered Enriques diagram with at most one satellite
vertex, then $\D$ is prime.
\end{Pro}
\begin{proof}
After \ref{roots}, we only need to prove 
that an ordered Enriques diagram $\D$ with one
satellite vertex, say $p_i(\D)$, is prime. If we are able to prove that
$\D^i$ is prime, the proof will be complete, because of
\ref{capsat}. So it is not restrictive to assume $i=r$.
Let $\D'\ne \D$ be such that $\Eff(\D')\subset \Eff(\D)$; 
we have to prove that $\dim \D'< \dim \D$. If 
$\dim \breve \D \ge \dim \breve \D' + 2$
then $\dim \D \ge \dim \D' + 1$ 
because $p_r(\D')$ is not a root (\ref{roots}), 
so we can assume $\dim \breve \D \le \dim \breve \D' + 1$.

Assume first $\dim \breve \D = \dim \breve \D' + 1$. As
$\breve \D$ has no satellite vertex and
$\dim \D= roots(\D)+frees(\D)$, 
\ref{roots} implies that there is an index $j$ satisfying
\begin{enumerate}
\item Either $p_j(\breve \D)$ is a root and
$p_j(\breve \D')$ is proximate to exactly one vertex or
$p_j(\breve \D)$ is proximate to exactly one vertex and
$p_j(\breve \D')$ is a satellite, and
\item For $k\ne j$, $p_k(\D')$ is proximate to $p_{\ell}(\D')$ if and only
if $p_k(\D)$ is proximate to $p_{\ell}(\D)$
\end{enumerate} 
Now, given a cluster $K\in \Cl(\breve \D')$,
if $p_k(\breve \D)$ is extremal among the vertices proximate to 
$p_j(\breve \D)$, $k>j$, then 
it easily follows by induction on $k-j$ that $E_j^{\D}$ 
and $E_k^{\D}$ meet transversely at a single point $p_{jk} \in S_K$.
So, if $\dim \breve \D = \dim \breve \D' + 1$ then
$p_r(\D')$ must be a satellite and we get 
$\dim \D = \dim \D' + 1$.

It only remains to consider the case $\dim \breve \D = \dim \breve \D'$.
By \ref{roots} it is clear that $\breve \D$ is prime, 
so $\breve \D'=\breve \D$. Let $p_j(\D)$ and
$p_k(\D)$ be the two vertices to which $p_r(\D)$ is
proximate. For every $K \in \Cl(\D')$,
$E^{\breve \D}_j$ and $E^{\breve \D}_k$
are irreducible in $S_{\breve K}$; therefore
$p_r(K)\in E^{\breve \D}_j \cap E^{\breve \D}_k$
is proximate to $p_j(K)$ and $p_k(K)$, and
$\D'=\D$, aganist the assumption.
\end{proof}

We introduce now invariants of pairs of diagrams $\D$, $\D'$
satisfying $\Eff(\D') \subset \Eff(\D)$, which will help
us to give new sufficient conditions for $\D$ to be prime.
Let $P$ and $P'$ be the proximity matrices of $\D$ and 
$\D'$ respectively, and let
$$M(\D,\D')=\left({P'}^{-1}P\right)^{T}{P'}^{-1}P \ .$$
Since $\Eff(\D') \subset \Eff(\D)$, all entries of
${P'}^{-1}P$ are greater or equal to zero; therefore
the $(j,k)$-th entry $m_{j,k}(\D,\D')$ of $M(\D,\D')$
is also greater or equal to zero. Moreover
it is easy to see, after \ref{compartir}, that
$m_{j,k}(\D,\D')>0$ if and only if for every cluster $K$ with
Enriques diagram $\D'$ the divisors $E_j^{\D}$ and
$E_k^{\D}$ share some component on $S_K$. For every
$i$, define
$$
\delta_i(\D,\D')=\# \left\{(k,j) \left| 
               \begin{matrix}  k >j \ge i, \, m_{j,k}(\D, \D') \ne 0 \\
                          p_k(\D) \text{ is extremal proximate to }
                          p_j(\D) 
               \end{matrix}
               \right. \right\}
$$
and 
$$
\epsilon_i(\D,\D')= (\dim \D - \dim \D') - (\dim \D^i -\dim {\D'}^i) 
           -\delta_i(\D, \D') \ .
$$ 

  A cluster $K$ (resp. an Enriques diagram $\Du$) is \emph{unibranched}
if it has only one root and no two points of $K$ (resp. vertices of
$\Du$) have the same immediate predecessor (i.e. there is at most one 
point of $K$ in the first neighbourhood of each point of $K$). It
is \emph{poli-unibranched} if it is the union of finitely many
unibranched clusters (resp. diagrams). Unibranched clusters
are also called bamboo clusters in the literature.
\begin{Rem}
\label{poliord}
  A unibranched cluster or Enriques diagram has only one
admissible ordering, in which each point or vertex
but the root is in the first neighbourhood of the
previous one.
A cluster or Enriques diagram which is not unibranched
has always more than one admissible ordering; for
poli-unibranched clusters and diagrams we shall (\ref{moltprimes})
use orderings mimicking the behaviour of unibranched
ones: an ordered poli-unibranched cluster or
Enriques diagram will be for us an ordered cluster or
Enriques diagram in which each point or vertex
but the roots is in the first neighbourhood or is a successor
of the previous one. 
\end{Rem}

\begin{Pro}
\label{moltprimes}
If $\D$ is a poli-unibranched ordered Enriques diagram
(i.e. for all $i$, $p_i(\D)$ is either a root or
proximate to $p_{i-1}(\D)$) then $\D$ is prime.
In particular, all unibranched Enriques diagrams are prime.
\end{Pro}

To prove Proposition \ref{moltprimes} we will need the following lemma:

\begin{Lem}
\label{escaleta}
Let $\D$ and $\D'$ be ordered Enriques diagrams with $r$ points
satisfying $\Eff(\D') \subset \Eff(\D)$, and let $i$ 
be an integer, $1 \le i \le r$.
Suppose that every point after $p_i(\D)$ is infinitely near
to $p_i(\D)$, is not proximate to any point before $p_i(\D)$
and is either satellite or proximate to $p_i(\D)$. Then
\begin{enumerate}
\item \label{delta} $\delta_i(\D, \D')=0$.
\item \label{epsilon} $\epsilon_i(\D,\D') \ge 0$, and 
$\epsilon_i(\D,\D')=0$ if and only if the proximity 
relations between points after $p_i$ are the same for $\D$ and $\D'$.
\item \label{irred} For every $K \in \Cl(\D')$, $E_i^{\D} \subset S_K$ is
irreducible.
\item \label{transvers} If $p_k(\D)$ is extremal proximate to 
$p_j(\D)$, $k>j \ge i$, then for every 
$K \in \Cl(\D')$, $E_j^{\D}$ 
and $E_k^{\D}$ meet transversely at a single point $p_{jk} \in S_K$.
\item \label{compta0} If $p_k(\D)$ is extremal proximate to 
$p_j(\D)$, $k>j > i$, then there is no other
pair $(j',k')$ such that $p_{k'}(\D)$ is extremal proximate to 
$p_{j'}(\D)$ and $E_{j'}^{\D} \cap E_{k'}^{\D} = p_{jk}$
in $S_K$, $K\in \Cl(\D')$.
\item \label{compta} For each $k>i$ such that
$p_k(\D)$ is extremal proximate to $p_i(\D)$
there are at most $\epsilon_i(\D,\D')$ 
indices $k' \ne k$ such that $p_{k'}(\D)$ is extremal proximate to 
$p_i(\D)$, and $E_i^{\D} \cap E_{k'}^{\D}=p_{ik}$
in $S_K$, $K\in \Cl(\D')$. 
\end{enumerate}
\end{Lem}

\begin{proof}
By induction on $r-i$. For $r=i$ there is nothing to prove.
For $r>i$, let $\breve \D$ and $\breve \D'$ be the diagrams 
obtained from $\D$ and $\D'$, respectively, by 
dropping the last point. Let $K$ be a cluster
in $\Cl(\D')$, and $\breve K \in \Cl(\breve \D')$ the cluster obtained
by dropping its last point. 
We distinguish two cases, depending on whether
$p_r(\D)$ is satellite or free and proximate to $p_i(\D)$.

If $p_r(\D)$ is satellite, then by equations (\ref{construccioXD})
$$p_r(K) \in E_j^{\breve \D} \cap E_k^{\breve \D}$$
for some $k > j \ge i$ such that $p_k(\breve \D)$ 
is extremal proximate to $p_j(\breve \D)$.
By part \ref{transvers} of the induction hypothesis 
we obtain that $\smash{E_j^{\breve \D}}$ and $\smash{E_k^{\breve \D}}$ 
meet transversely in a single point; this point must be
$p_r(K)$, which is therefore a satellite, so
$p_r(\D')$ is a satellite too.
Denoting by $\pi$ the blowing up of $p_r(K)$, 
the virtual exceptional divisors are 
\begin{align*}
E_r^{\D}&= E_r = \pi^{*}(p_r(K)), \\
E_j^{\D}&= \pi^{*}(E_j^{\breve \D}) -E_r, \\
E_k^{\D}&= \pi^{*}(E_k^{\breve \D}) -E_r, \\
E_n^{\D}&= \pi^{*}(E_n^{\breve \D}) \text{ if } n \ne j, k, r. 
\end{align*}
So parts \ref{irred}, \ref{transvers}, \ref{compta0}
and \ref{compta} of the claim are consequences of
the corresponding parts of the induction hypothesis.
Part \ref{delta} is immediate from \ref{transvers}
and \ref{epsilon} follows from \ref{compta}. Note
that in this case 
$\epsilon_i(\D,\D')=\epsilon_i(\breve \D, \breve \D')$.

If $p_r(\D)$ is free and proximate to $p_i(\D)$,
then we have two possibilites: either $p_r(K)$ is
also free, in which case 
$\epsilon_i(\D,\D')=\epsilon_i(\breve \D, \breve \D')$ and
all claims follow immediately from the induction hypothesis,
or $p_r(K)= E_i^{\breve \D} \cap E_k^{\breve \D}$, in
which case more interesting phenomena occur.
Denoting by $\pi$ the blowing up of $p_r(K)$,
the virtual exceptional divisors are 
\begin{align*}
E_r^{\D}&= E_r = \pi^{*}(p_r(K)), \\
E_i^{\D}&= \pi^{*}(E_i^{\breve \D}) -E_r, \\
E_n^{\D}&= \pi^{*}(E_n^{\breve \D}) \text{ if } n \ne i, r.
\end{align*}
This, together with the induction hypothesis, proves parts 
\ref{irred} and \ref{compta0} of the claim, but it also means  
that $E_r$ is a component of every divisor $E_n^{\D}$ such that 
$p_r(K)= E_i^{\breve \D} \cap E_n^{\breve \D}$.
By part \ref{compta} of the induction hypothesis
we know that there can be at most $\epsilon_i(\breve \D,\breve \D')+1$
such divisors, including $n=k$, and for every one of them
we have $E_i^{\D} \cap E_n^{\D}= E_i^{\D} \cap E_r^{\D}$.
$p_r(\D)$ is not proximate to any $p_n(\D)$, so
\ref{transvers} has been proved also, and therefore
\ref{delta}. As $\delta_i(\D,\D')=0$, the fact that
$p_r(K)$ and $p_r(\D')$ are satellite whereas $p_r(\D)$
is free means that 
$\epsilon_i(\D,\D')=\epsilon_i(\breve \D,\breve \D')+1$
thus proving parts \ref{compta} and \ref{epsilon}.
\end{proof}

\begin{proof}[Proof of \ref{moltprimes}]
Let $\D_1 \ne \D$ be an ordered Enriques diagram with 
$\Eff(\D_1) \subset \Eff(\D)$.
We have to prove that $\dim \D_1 < \dim \D$. Fix a cluster $K \in \Cl(\D_1)$.
For every $i=1, 2, \dots, r$
let $r_i$ be such that $p_{r_i}(\D)$ is the last point in $\D$ proximate
to $p_i$, and consider the Enriques diagram $\D^{r_i}$ and
the cluster $K^{r_i}$. 

The hypothesis that $\D$ is unibranched or poli-unibranched 
says that for every $i$ and for every satellite point
$p_j(\D)$, proximate to $p_i(\D)$, all vertices
between $p_i(\D)$ and $p_j(\D)$ are
infinitely near to $p_i(\D)$, and satellite or 
proximate to $p_i(\D)$.
This, toghether with lemma \ref{escaleta}, 
show that $\epsilon_i(\D^{r_i},\D_1^{r_i}) \ge 0$ 
for all $i$.
On the other hand, if $i$ is the least integer such that $\D$
and $\D_1$ differ at the point $p_i$ then 
$\epsilon_{i-1}(\D^{i},\D_1^{i}) > 0$, which clearly implies
that $\dim \D^i - \dim \D_1^i >0$. Now we can prove
that $\dim \D - \dim  \D_1 >0$ by induction on $r-i$.
Assume we have proved $\dim \D^j - \dim \D_1^j >0$
for all $j$, $i \ge j <r$, and let us prove it for $j=r$.
There are two possibilities, either $p_r(\D)$ is a root, or
there is $j<r$ with $r=r_j$. If $p_r(\D)$ is a root, then
$\dim \D - \dim \D_1 \ge \dim \breve \D - \dim \breve \D_1 >0$,
by the induction hypothesis.
If $r=r_j$ and $i>j$ then by \ref{escaleta} 
$\epsilon_j(\D,\D_1)>0$ which implies $\dim \D - \dim  \D_1 >0$
and we are done. If $r=r_j$ and $i \le j$ then 
$\epsilon_j(\D,\D_1) \ge 0$, and the induction hypothesis says
$\dim \D^j - \dim \D_1^j >0$; the two inequalities together give
$$
\dim \D - \dim  \D_1 \ge \dim \D^j - \dim \D_1^j + \epsilon_j(\D,\D_1) >0
$$
thus finishing the proof.
\end{proof}


\vskip 2mm

As we see, there are many interesting prime Enriques diagrams $\D$, and 
for them we know exactly which diagrams are specializations of $\D$.
If $\D$ is a non-prime diagram, or if we do not know whether it
is prime or not, we do not have necessary and sufficient 
conditions for specialization, but we can still give some information,
in the form of a sufficient condition.

Given a map 
$\sigma: \{1, 2, \dots, r\} \longrightarrow \{1, 2, \dots, r \}$
we define the matrix of $\sigma$ to be $\Sigma=(\sigma_{ij})$ where
$$
\sigma_{ij}=
\begin{cases}
1 & \text{ if } \sigma(j)=i \\
0 & \text{ otherwise.}
\end{cases}
$$

\begin{Pro}
\label{esten}
Let $\D_1$, $\D_2$ be two Enriques diagrams with $r$ vertices,
$s$ an integer,  $1 \le s \le r$, and 
$\sigma: \{1, 2, \dots, r\} \rightarrow \{1, 2, \dots, r\}$
a map such that $\sigma(i)=i \, \forall i>s$. 
Let $P_1$ and $P_2$ be the proximity matrices
of $\D_1$ and $\D_2$ respectively, and $\Sigma$ the matrix
of $\sigma$.
Consider the diagrams
$\D_1^s, \D_2^s$ of the first $s$ points of $\D_1$ and $\D_2$ respectively,
and assume
\begin{enumerate}
   \item \label{partida} $\D_1^s \rightsquigarrow \D_2^s$,
   \item \label{matrius} the matrix $P_2^{-1}P_1-\Sigma$ has no negative entries,
   \item \label{nosalta} if $p_i(\D_1)$, $i>s$ is a satellite proximate
         to $p_j(\D_1)$ and $p_k(\D_1)$, then $j,k>s$,
   \item \label{constancia} for every $j>s$,  
         $p_i(\D_1) \prec p_j(\D_1)$ if and only if
         $p_{\sigma(i)}(\D_2) \prec p_j(\D_2)$.
\end{enumerate}
Then $\D_1 \rightsquigarrow \D_2$. Moreover, for every $i>s$
and $K\in \Cl(\D_2)$, the virtual divisor $E_i^{\D_1} \subset S_K$
is irreducible.
\end{Pro}

\begin{Exa}
Consider the Enriques diagrams with eight vertices of
figure \ref{doscusps2}.
\begin{figure}
  \begin{center}
    \mbox{\includegraphics{esten.eps}}
    \caption{$\D_1 \rightsquigarrow \D_2$.}
    \label{doscusps2}
  \end{center}
\end{figure}
It is clear, after example \ref{ex1}, that $\D$ is not
prime. However, we can prove that
$\D_1 \rightsquigarrow \D_2$, by using Proposition \ref{esten}
with $s=2$, and $\sigma(1)=\sigma(2)=2$. $\D_1^2$ has no
satellite vertices, so it is prime and 
$\D_1^2 \rightsquigarrow \D_2^2$; one can then easily check
the other conditions needed to apply \ref{esten}.
\end{Exa}

\begin{proof}
Let $K\in \Cl(\D_2)$; we have to see that $K \in \overline{\Cl(\D_1)}$.
We shall do it by induction on $r-s$. For $r=s$ the claim follows
immediately from hypothesis \ref{partida}.
For $r>s$ consider the Enriques diagrams $\breve \D_1$
and $\breve \D_2$ and the cluster $\breve K$ obtained by 
dropping the last vertex of $\D_1$, $\D_2$ and the
last point of $K$ respectively.
The induction hypothesis is that 
$$\breve K\in \overline{\Cl(\breve \D_1)}$$ and for every $i>s$
and $K\in \Cl(\breve \D_2)$, $E_i^{\breve \D_1}$ is irreducible. 

We shall now fix our attention to the family 
$$                       
\begin{CD}
    S_{\overline{\Cl(\breve \D_1)}} @>\psi_{r-1}>> \overline{\Cl(\breve \D_1)}.
\end{CD}
$$
To lighten notations a little bit, let $T=\overline{\Cl(\breve \D_1)}$.
Since $\breve K = \psi_{r-1}(K)$ (\ref{clustersfam}) the induction 
hypothesis implies $K \in S_T$.
On the other hand, $\psi_{r-1}$ is smooth and has irreducible fibers, and
$T$ is irreducible, so
$S_T$ is irreducible.
If $p_r(\D_1)$ is a root, then $\Cl(\D_1)$ can be identified with
an open subset of $S_T$, and due to the 
irreducibility of $S_T$ we get
$\overline{\Cl(\D_1)}=S_T$ and $K \in \overline{\Cl(\D_1)}$.

The divisors $F_{i,T}^{\breve \D_1}$ 
are relative divisors over $T$ 
so every component dominates $T$ 
(because $T$ is irreducible). 
Also, for every cluster $K' \in \Cl(\D_1) \subset S_T$  
$E_i^{\breve \D_1}=F_{i,T}^{\breve \D_1}|_{S_K}$ are irreducible,
therefore the $F_{i,T}^{\breve \D_1}$ are irreducible.
If $p_r(\D_1)$ is free and proximate to $p_i(\D_1)$, 
then $\Cl(\D_1)$ is
an open subset of $F_{i,T}^{\breve \D_1}$.
Due to the irreducibility of 
$F_{i,T}^{\breve \D_1}$ this implies
$\overline{\Cl(\D_1)}=F_{i,T}^{\breve \D_1}$.
On the other hand, by hypothesis \ref{constancia}
$p_{\sigma(i)}(K) \prec p_r(K)$, so  
$$
K \in F_{\sigma(i),T}^{\breve \D'}\ ,
$$
and by hypothesis \ref{matrius}
$$
F_{\sigma(i),T}^{\breve \D_2} \subset
F_{\sigma(i),T}^{\breve \D_1}\ ,
$$
so we obtain $K \in F_{\sigma(i),T}^{\breve \D_1}$,
as wanted.

Suppose now that $p_r(\D_1)$ is a satellite, and
$p_i(\D_1)$, $p_j(\D_1)$ are the two vertices to which
$p_r(\D_1)$ is proximate.
$\overline{\Cl(\D_1)} \rightarrow T$ is
onto, so there is at least one cluster $K' \in \overline{\Cl(\D_1)}$
with $\psi_{r-2}(K')=\breve K$. 
Because of \ref{constancia}, $p_{\sigma(i)}(K) \prec p_{r}(K)$
and $p_{\sigma(j)}(K) \prec p_{r}(K)$, so
$$K\in E_{\sigma(i)}^{\breve \D_2} \cap E_{\sigma(j)}^{\breve \D_2}$$
and by hypothesis \ref{matrius}
$$
E_{\sigma(i)}^{\breve \D_2} \cap E_{\sigma(j)}^{\breve \D_2} \subset
E_{\sigma(i)}^{\breve \D_1} \cap E_{\sigma(j)}^{\breve \D_1} \ .
$$
On the other hand, clearly $K'\in \Eff(\D_1)$,
so we get 
$$ p_r(K), p_r(K') \in E_i^{\breve \D_1} \cap E_j^{\breve \D_1} \ .$$
But $E_i^{\breve \D_1}$ and $E_j^{\breve \D_1}$ are irreducible,
because of hypothesis \ref{nosalta} and the induction hypothesis, 
so its intersection is the single point
$p_r(K)=p_r(K')=E_i^{\breve \D_1} \cap E_j^{\breve \D_1}$,
and $K=K' \in \overline{\Cl(\D_1)}$.

The fact that  the only virtual divisor 
$E_j^{\breve \D_1}$ in $S_{\breve K}$
which has $\tilde E_i$ as a component is $E_i^{\breve \D_1}$ is 
immediate, from the induction
hypothesis and the fact that $p_r(K)$ belongs to $E_i^{\breve \D_1}$
if and only if $p_r(\D_1)$ is proximate to $p_i(\D_1)$.
\end{proof}

\begin{Pro}
\label{esten2}
Let $\D_1$, $\D_2$ be two Enriques diagrams with $r$ vertices,
$s$ an integer,  $1 \le s \le r$, and 
$\sigma: \{1, 2, \dots, r\} \rightarrow \{1, 2, \dots, r\}$
a map such that $\sigma(i)=i \, \forall i>s$.
Consider the diagrams
$\D_1^s, \D_2^s$ of $s$ points obtained by dropping the
last $r-s$ points of $\D_1$ and $\D_2$ respectively, and assume
\begin{enumerate}
   \item $\D_1^s \rightsquigarrow \D_2^s$,
   \item $\Eff(\D_2) \subset \Eff(\D_1)$,
   \item $E_i^{\D_1^s} - \tilde E_{\sigma(i)}$ is effective in $S_K$
         for all $K \in \Cl(\D_2^s)$,
   \item There is no satellite point $p_i(\D_1)$ with
         $i>s$.
\end{enumerate}
Then $\D_1 \rightsquigarrow \D_2$.
\end{Pro}
\begin{proof}
The same as for Proposition \ref{esten}; the divisors
$E_i^{\breve \D_1}$ can have common components but this has
no consequence since there appear no satellite points.
\end{proof}

\subsection{Extension of clusters}

Extension is a particular kind of specialization of clusters
which was first introduced by Greuel, Lossen and Shustin
in \cite{GLS98}, and is basic in their construction of singular
curves of low (asimptotically proper) degree. In their approach
it is in fact a specialization of zero-dimensional schemes,
presented using a specialization of singular curves. They
give also a description of extension in terms of a specialization
of \emph{weighted consistent} clusters. As we will shortly see,
the specialization of clusters can be justified independently of
the weights and with no reference to any particular curve, and
by doing so, slightly generalized.

Let $\Du$ be an unordered Enriques diagram with $r$ vertices, 
$q \in V(\Du)$ one of its free vertices, proximate to its predecessor
$\bar q$ only. The \emph{extension} of $\Du$ at $q$ is the
Enriques diagram $\Du(q)$ which has:
\begin{itemize}
   \item The same vertices as $\Du$, plus a new free vertex $q'$
         inserted between $\bar q$ and $q$. That is,
         the predecessor of $q'$ is $\bar q$ and the predecessor
         of $q$ is $q'$.
   \item $q'\prec p \text{ in }\Du(q)  \Longleftrightarrow 
         q \le p \text{ and } \bar q \prec p \text{ in }\Du$.
   \item $\bar q\prec p \text{ in }\Du(q)  \Longleftrightarrow 
         q \not \le p \text{ and } \bar q \prec p \text{ in }\Du$.
   \item All proximities $p \prec p'$, $p \ne \bar q, q'$ remain unchanged.
\end{itemize}
We consider also the Enriques diagram $\Du_\bullet$ which has the
same vertices and proximities as $\Du$ plus one root. 
The specialization introduced in \cite{GLS98} (in which
all vertices preceding $q$ are assumed free) gives a 
flat family of zero-dimensional schemes, whose general member
is the scheme of a cluster whose Enriques diagram is
$\Du_\bullet$ and has a special member which is the scheme of 
a cluster whose Enriques diagram is $\Du(q)$.

\begin{Pro}
Choose an admissible ordering $\D$ for $\Du$, such that
all vertices after $q$ in this ordering are infinitely
near to $q$. Let $\D_\bullet$ and $\D(q)$ be the orderings 
induced by this ordering on $\Du_\bullet$ and $\Du(q)$ 
with the additional condition that the new root (resp. $q'$) is
the vertex immediately preceding $q$.
Then $\D_\bullet \rightsquigarrow \D(q)$. 
\end{Pro}
\begin{proof}
Suppose that $q$ is the $s$-th vertex of $\D$,
and $\bar q$ is the $t$-th. 
Then $q'$ is the $s$-th vertex of $\D(q)$,
and the diagrams $\D_\bullet^s$, $\D(q)^s$ differ only in their
last point, which is a root for $\D_\bullet^s$.
Let $T$ be the irreducible variety
$$
T=\overline{\Cl(\D^{s-1})}=\overline{\Cl(\D_\bullet^{s-1})}
 =\overline{\Cl(\D(q)^{s-1})}
$$
as in the proof of \ref{esten}. $p_s(\D_\bullet^s)$
is a root, and therefore $\Cl(\D_\bullet^s)=S_T$, which
implies that every cluster of $s$ points whose
$s-1$ first points have Enriques diagram $\D^s$
belongs to $\Cl(\D_\bullet^s)$. In particular, this 
applies to clusters in $\Cl(\D(q)^s)$, so we 
get $\D_\bullet^s \rightsquigarrow \D(q)^s$. Let
$\sigma: \{1, 2, \dots, r+1\} \longrightarrow \{1, 2, \dots, r+1 \}$
be defined by $\sigma(t)=s$ and $\sigma(i)=i \, \forall i \ne t$.
Now the proof follows as in Proposition \ref{esten},
proving by induction on $r-s$ that $\D_\bullet \rightsquigarrow \D(q)$,
and that for every $i \ge s$ and $K\in \Cl(\D(q))$, 
the virtual divisor $E_i^{\D_\bullet} \subset S_K$
is irreducible, and for every $k$ such
that $p_k(\D)$ is maximal proximate to $\bar q$
$$
E_k^{\D_\bullet} \cap E_t^{\D_\bullet} = \tilde E_k \cap \tilde E_s .
$$
\end{proof}

\section{Enriques diagrams do not stratify $X_{\MakeLowercase{r}-1}$}
\label{nostrata}

In this section we prove that the subvarieties
$\Cl(\D)$ do not constitute a stratification
on $X_{r-1}$, because it is not true that  
$\overline{\Cl(\D_1)} \cap \Cl(\D_2) \ne \emptyset$ 
implies $\Cl(\D_2) \subset \overline{\Cl(\D_1)}$.
To prove it, we shall show an explicit example
of two diagrams $\D_1$ and $\D_2$ with
$\overline{\Cl(\D_1)} \cap \Cl(\D_2) \ne \emptyset$ 
and $\Cl(\D_2) \not\subset \overline{\Cl(\D_1)}$.

Consider the Enriques diagrams with seven points of
figure \ref{sisfigs}.
It is clear, after example \ref{ex2}, that $\D_1$ is not
prime, so $\Eff(\D_1)$ has at least two components.
Moreover, it is also clear that 
$\Cl(\D_2)\subset\Eff(\D_2)\subset \Eff(\D_1)$,
and $\dim \D_1=\dim \D_2=7$, so 
$\Cl(\D_2)\not\subset \overline{\Cl(\D_1)}$. We shall
see that $\Cl(\D_2) \cap \overline{\Cl(\D_1)} \ne \emptyset$,
To simplify matters, we work on $\P^2$, and we shall
prove a stronger statement:

\begin{Pro}
\label{pronostrata}
Let $S=\P^2$, and let $K$ be
a cluster whose Enriques diagram is $\D_K$
(see figure \ref{sisfigs}).
  \begin{enumerate}
  \item $\overline{\Cl(\D_2)}$ is a component of $\Eff(\D_1)$.
  \item \label{cent} $K\in \overline{\Cl(\D_2)} \cap \overline{\Cl(\D_1)}$.
  \item \label{by} $\D_1 \rightsquigarrow \D_K$ and $\D_2 \rightsquigarrow \D_K$.
  \item $\overline{\Cl(\D_2)} \cap \overline{\Cl(\D_1)}$ has five components going
through $K$ (one of them double) which are $\Eff(\D_{x_4})$,
$\Eff(\D_{x_6})$, $\Eff(\D_{y_6})$ (see figure \ref{sisfigs})
and two components contained in $\Cl(\D_2)$.
  \end{enumerate}
\end{Pro}

\begin{proof}
Note that \ref{by} is an immediate consequence of
\ref{cent}.
The method of proof uses coordinates and computations
with the computer algebra package {\bf Macaulay}.

Let $V_1 \cong k^2$ be an affine neigbourhood of
$p_1(K)$ in $\P^2$ with coordinates $(x_1, y_1)$ such
that $p_2(K)$ lies on the $x_1$ axis, $y_1=0$.

After the first blowing up, there is a neighbourhood
of $p_2(K)\in X_1$ isomorphic to $k^4$, where the
coordinate functions are $x_1 \circ \pi_1$, $y_1 \circ \pi_1$,
$x_2$ and $y_2$. To simplify
notation, we write $x_1 =x_1 \circ \pi_1$
and $y_1 =y_1 \circ \pi_1$.
It is not hard to see that $x_2$, $y_2$
can be chosen in such a way that
the restricted blowing-up $\pi_1:V_2 \rightarrow V_1$ 
is given by
$$
\pi_1(x_1, y_1, x_2, y_2)=(x_1+x_2, y_1+x_2 y_2) ,
$$
and the equation of $F_1$ and of $\Cl(\D_2^2)$ in this
neighbourhood is $x_2=0$.
 
After each blowing up, there is an affine neigbourhood
of $p_i(K)$ in which we have coordinates and can keep
track of the equations of the subvarieties we are 
interested in.
The iterated blowing-up process continues until we reach $p_7(K)$.
 We summarize the process in the table \ref{taula},
denoting by $F_i(\D)$ the effective divisor on
$S_\Eff(\D)$ which is the restriction of $F_i^{\D}$.

\begin{figure}
  \begin{center}
    \mbox{\includegraphics{nostrat1.eps}}

\vskip 2cm

    \mbox{\includegraphics{nostrat2.eps}}

\vskip 2cm

    \mbox{\includegraphics{nostrat3.eps}}
    \caption{$\Cl(\D_2)\cap \overline{\Cl(\D_1)}\ne \emptyset$,
             but $\D_1 \not\rightsquigarrow \D_2$.}
    \label{sisfigs}
  \end{center}
\end{figure}

\begin{table}
\label{taula}
\small
\begin{tabular}{|p{\taulength}|}
\hline
$\pi_1(x_1, y_1, x_2, y_2)=(x_1+x_2, y_1+x_2 y_2)$

\medskip

$I_V(F_1)=(x_2)\, ,$

$I_V(\Eff(\D_2^2))=(x_2)$ \\
\hline
$\pi_2(x_1, y_1, \dots, x_3, y_3)=(x_1, y_1, x_2+x_3 y_3 , y_2+ y_3)$ 

\medskip

$I_V(F_2)=(y_3)\, , $

$I_V(F_1)=(x_2 \circ \pi_1)=(x_2 + x_3 y_3)\, , $

$I(F_1(\D_1^2))=(x_2 + x_3 y_3)\, , $

$I(F_1(\D_2^2))=(x_3)\, ,$

\medskip

$I_V(\Eff(\D_1^3))=(y_3)\, ,$

$I_V(\Eff(\D_2^3))=(x_2, x_3, y_3)$ \\
\hline
$\pi_3(x_1, y_1, \dots, x_4, y_4)=(x_1, y_1, x_2, y_2 , x_3+x_4 y_4 , y_3+ y_4)$ 

\medskip

$I_V(F_3)=(y_4)\, , $

$I(F_1(\D_1^3))=(x_2 + (x_3 + x_4 y_4) y_4)\, ,$

$I(F_1(\D_2^3))=(x_4)\, ,$

\medskip

$I_V(\Eff(\D_1^4))=(y_3,x_2+(x_3+x_4 y_4)y_4)\, ,$

$I_V(\Eff(\D_2^4))=(x_2, x_3, y_3, y_4)$ \\
\hline
$\pi_4(x_1, y_1, \dots, x_5, y_5)=
(x_1, y_1, \dots, x_3, y_3 , x_4+x_5 , y_4+ x_5 y_5)$ 

\medskip

$I_V(F_4)=(x_5)\, , $

$I(F_3(\D_2^4))=(y_5)\, ,$

$I(F_1(\D_1^4))=(y_4^2 + x_3 y_5 + 2 x_4 y_4 y_5 + 2 x_5 y_4
y_5 + x_4 x_5 y_5^2 + x_5^2 y_5^2)\, ,$

\medskip

$I_V(\Eff(\D_1^5))=(y_3, x_2 + (x_3+x_4 y_4)y_4, x_5,
  y_4^2 + x_3 y_5 + 2 x_4 y_4 y_5)\, ,$

$I_V(\Eff(\D_2^5))=(x_2, x_3, y_3, y_4, x_5)$ \\
\hline
$\pi_5(x_1, y_1, \dots, x_6, y_6)=
(x_1, y_1, \dots, x_4, y_4 , x_5+x_6 , y_5+ x_6 y_6)$ 

\medskip

$I_V(F_5)=(x_6)\, , $

$I(F_3(\D_2^5))=(y_5+x_6 y_6)\, ,$

$I(F_1(\D_1^5))=(2 y_4 y_5 + x_4 y_5^2 + x_6 y_5^2 + x_3 y_6 + 
 2 x_4 y_4 y_6 + 2 x_6 y_4 y_6 + 2 x_4 x_6 y_5 y_6 + 
  2 x_6^2 y_5 y_6 + x_4 x_6^2 y_6^2 + x_6^3 y_6^2)\, ,$

\medskip

$I_V(\Eff(\D_1^6))=(y_3, x_2 + (x_3+x_4 y_4) y_4, x_5,
  y_4^2 + x_3 y_5 + 2 x_4 y_4 y_5, 2 y_4 y_5 + x_4 y_5^2 + 
  x_6 y_5^2 + x_3 y_6 + 2 x_4 y_4 y_6 + 2 x_6 y_4 y_6 + 2 x_4 x_6 y_5 y_6 + 
  2 x_6^2 y_5 y_6 + x_4 x_6^2 y_6^2 + x_6^3 y_6^2)\, ,$

$I_V(\Eff(\D_2^6))=(x_2, x_3, y_3, y_4, x_5, y_5 + x_6 y_6)$ \\
\hline
$\pi_6(x_1, y_1, \dots, x_7, y_7)=
(x_1, y_1, \dots, x_5, y_5 , x_6+x_7 , y_6+ x_7 y_7)$ 

\medskip

$I_V(F_6)=(x_7)\, , $

$I(F_1(\D_1^6))=(y_5^2 + 2 y_4 y_6 + 2 x_4 y_5 y_6 + 
  4 x_6 y_5 y_6 + 2 x_7 y_5 y_6 + 2 x_4 x_6 y_6^2 + 3 x_6^2 y_6^2 + 
  x_4 x_7 y_6^2 + 3 x_6 x_7 y_6^2 + x_7^2 y_6^2 + x_3 y_7 + 
  2 x_4 y_4 y_7 + 2 x_6 y_4 y_7 + 2 x_7 y_4 y_7 + 2 x_4 x_6 y_5 y_7 + 
  2 x_6^2 y_5 y_7 + 2 x_4 x_7 y_5 y_7 + 4 x_6 x_7 y_5 y_7 + 
  2 x_7^2 y_5 y_7 + 2 x_4 x_6^2 y_6 y_7 + 2 x_6^3 y_6 y_7 + 
  4 x_4 x_6 x_7 y_6 y_7 + 6 x_6^2 x_7 y_6 y_7 + 2 x_4 x_7^2 y_6 y_7 + 
  6 x_6 x_7^2 y_6 y_7 + 2 x_7^3 y_6 y_7 + x_4 x_6^2 x_7 y_7^2 + 
  x_6^3 x_7 y_7^2 + 2 x_4 x_6 x_7^2 y_7^2 + 3 x_6^2 x_7^2 y_7^2 + 
  x_4 x_7^3 y_7^2 + 3 x_6 x_7^3 y_7^2 + x_7^4 y_7^2)\, ,$

\medskip

$I_V(\Eff(\D_1^7))=(y_3, x_2 + (x_3+x_4 y_4) y_4, x_5,
  y_4^2 + x_3 y_5 + 2 x_4 y_4 y_5, 2 y_4 y_5 + x_4 y_5^2 + 
  x_6 y_5^2 + x_3 y_6 + 2 x_4 y_4 y_6 + 2 x_6 y_4 y_6 + 2 x_4 x_6 y_5 y_6 + 
  2 x_6^2 y_5 y_6 + x_4 x_6^2 y_6^2 + x_6^3 y_6^2,   y_5^2 + 
  2 y_4 y_6 + 2 x_4 y_5 y_6 + 4 x_6 y_5 y_6 + 2 x_4 x_6 y_6^2 + 
  3 x_6^2 y_6^2 + x_3 y_7 + 2 x_4 y_4 y_7 + 2 x_6 y_4 y_7 + 
  2 x_4 x_6 y_5 y_7 + 2 x_6^2 y_5 y_7 + 2 x_4 x_6^2 y_6 y_7 + 
  2 x_6^3 y_6 y_7)\, ,$

$I_V(\Eff(\D_2^7))=(x_2, x_3, y_3, y_4, x_5, y_5 + x_6 y_6, x_7)$ \\
\hline
\end{tabular}
\caption{Equations of $\Eff(\D_1^i)$, $\Eff(\D_2^i)$.}
\end{table}

It is worth noting that  
$\overline{\Cl(\D_1^5)}$ is singular. This fact 
seems to be related to the presence of non reduced
virtual exceptional divisors (see example \ref{ex2}). 

Once we get the equations in $V_7$ of $\Eff(\D_1)$ and
$\Eff(\D_2)=\overline{\Cl(\D_2)}$ we use Macaulay
to obtain equations in $V_7$ of the residual component
$\Eff(\D_1) \setminus \overline{\Cl(\D_2)}= \overline{\Cl(\D_1)}$.
They turn out to be 
\begin{gather*}
I_V(\overline{\Cl(\D_1)})=(y_3, x_5, x_7,
   y_4^2+ x_3y_5 + 2 x_4y_4y_5, 
   x_4y_4^2+x_3y_4+x_2, \\
   x_4x_6^2y_6^2 + x_6^3y_6^2 + 2x_4x_6y_5y_6 + 2x_6^2y_5y_6 + 
   x_4y_5^2 + x_6y_5^2  + x_4y_4y_6 + 2x_6y_4y_6 + \\ + 2y_4y_5 + x_3y_6, 
   2x_4x_6^2y_6y_7 + 2x_6^3y_6y_7 + 2x_4x_6y_6^2 + 3x_6^2y_6^2 + 
   2x_4x_6y_5y_7 + \\ + 2x_6^2y_5y_7 + 2x_4y_5y_6 + 4x_6y_5y_6 + 2x_4y_4y_7 +
   2x_6y_4y_7 + y_5^2 + 2y_4y_6 + x_3y_7, \\ 
   4x_4^2x_6y_5y_6y_7 - 4x_6^3y_5y_6y_7 + 4x_4^2y_5y_6^2 + 
   2x_4x_6y_5y_6^2 - 8x_6^2y_5y_6^2 - 4x_4x_6y_5^2y_7 - \\
 - 4x_6^2y_5^2y_7 + 2x_4x_6y_4y_6y_7 + 2x_6^2y_4y_6y_7 -
   2x_4y_5^2y_6 - 12x_6y_5^2y_6 + 4x_4y_4y_6^2 + \\
 + 6x_6y_4y_6^2 - 2x_6y_4y_5y_7 - 4y_5^3 +
   2y_4y_5y_6 + y_4^2y_7 + x_3y_5y_7) \, .
\end{gather*}
In particular we see that $K$ (which has all its
coordinates equal to zero) belongs to this variety,
so we have $K\in \overline{\Cl(\D_2)} \cap \overline{\Cl(\D_1)}$,
as claimed above. If we now intersect the closures
of $\Cl(\D_1)$ and $\Cl(\D_2)$ the result is
\begin{gather*}
I_V(\overline{\Cl(\D_1)} \cap \overline{\Cl(\D_2)})=
(x_2, x_3, y_4, y_5 + x_6y_6, 
 x_4x_6y_6^2(x_4+y_6)(y_6+ x_6y_7)) \, .
\end{gather*}
This intersection has five components going
through $K$, one of them double, which correspond
to the five factors of the last generator of
$I_V(\overline{\Cl(\D_1)}$. The first three
components can be identified as the varieties of some
Enriques diagram because table \ref{taula} 
shows that the vanishing of $x_4$,
$x_6$ and $y_6$ can be interpreted as proximity
relations satisfied by points of the clusters.
Explicitly, these subvarieties are $\Eff(\D_{x_4})$,
$\Eff(\D_{x_6})$, and $\Eff(\D_{y_6})$. The other two
components intersect $\Cl(\D_2)$ nontrivially. The
component corresponding to the factor $x_4+x_6$
contains clusters whose Enriques diagram is $\D_2$;
the equation can be interpreted as imposing to
the fourth and sixth points (which lie on
$\tilde E_3$) a special position with
respect to the two satellite points which also
lie on $\tilde E_3$.                          
The last component dominates $\Cl(\breve \D_2) \cap V_6$
and contains clusters whose Enriques diagram is $\D_2$
also. Indeed, the coordinates in $V_6$ of a cluster 
$\breve K' \in \Cl(\breve \D_2) \cap V_6$ have
$x_6, y_6 \ne 0$, so choosing $p_r(K') \in S_{\breve K}$
with coordinates $x_7=0$, $y_7=y_6/x_6$ we
obtain a cluster in $\overline{\Cl(\D_1)}$
whose Enriques diagram is $\D_2$.

\end{proof}

\section{A functorial approach}
\label{functorial}

The iterated blowing-ups $X_r$ can be studied from
a functorial point of view following the ideas of
Harbourne in \cite{Har88}. Let $S$ as before be
a fixed smooth projective surface. The objects of
interest are now the surfaces obtained by blowing
up $r$ proper or infinitely near points on $S$, that is,
the surfaces $S_K$ where $K$ has $r$ points. One looks
for a moduli space for those surfaces, that is,
a parametrization with one point for each isomorphism
class. Unfortunately, such a space does not exist, 
but there is a natural parametrization
of these varieties which has many of the functorial
properties a moduli space would have, and moreover
it contains information on automorphisms and jumping
of structure which a genuine moduli space would not
have. This parametrization is given by 
$\psi_r: X_r \rightarrow X_{r-1}$. In this section
we shall obtain some consequences of
\ref{efclo} and \ref{loclo} by looking at the
functorial properties of $X_r$ and its subvarieties.


Let $\F_r$ be the contravariant functor from the 
category $\Sch$ of schemes of finite type over
$k$ to the category $\Sets$, where for any such
scheme $T$, $\F_r(T)$ is the set of families of
surfaces $S_K$ parametrized by $T$, i.e., 
the set of $T$--isomorphism
classes of schemes proper and smooth over $T$
whose fibres at closed points are surfaces
obtained by blowing up $S$ at $r$ points.
We consider also the contravariant functor 
$\B_r$ from $\Sch$ to $\Sets$ of ordered 
blowing-ups of $r$ points of $S$. An ordered
blowing-up of $r$ $T$--points of $S$ is a sequence
of morphisms of $T$--schemes 
$\pi_i:Y_i \rightarrow Y_{i-1}$ where $Y_0=S \times T$.
Then $\B_r(T)$ is defined to be the set of ordered blowing-ups
of $r$ $T$--points of $S$.
Remark that given an ordered blowing-up of $r$ $T$--points of $S$
for each closed $t \in T$ the fiber $(Y_r)_t$ is a surface $S_K$
where $K$ is a well defined ordered cluster determined by the order
of the blowing-ups. We will denote this cluster $K(t)$.

Following \cite{Har88}, we say that
a family $\psi: {\mathcal S} \rightarrow T$ together
with a couple of morphisms 
$s, t : {\mathcal I} \rightarrow T$
\emph{solves the moduli problem} for the functor
$\F$, provided that the following conditions hold:
\begin{itemize}
\item The family $\psi$ is
versal and any family locally for the \'etale topology
comes from it; i.e., if $(X \rightarrow Y) \in \F(Y)$
for some scheme $Y$, then each point $y$ of $Y$ has
an \'etale neighborhood $V$ and a morphism $V \rightarrow T$
such that $X_V \rightarrow V$ comes from $\psi$
by base change $V \rightarrow T$.
\item The scheme $\mathcal I$ represents the functor
$\Isom$ where for any shcheme $U$ and morphisms
$f_i: U \rightarrow T$, $i=1,2$, $\Isom (U)$ is the
set of $U$--isomorphisms from $f_1^{*}(\psi)$ to
$f_2^{*}(\psi)$. Each closed point $i$ of $I$ corresponds
to an isomorphism between two surfaces of the
family ${\mathcal S} \rightarrow T$, namely the 
ones parametrized by $s(i)$ (source) and $t(i)$
(target).
\end{itemize}

Harbourne proved (\cite[I.2]{Har88}) that the 
variety $X_{r-1}$ represents the functor $\B_r$,
the universal family being 
$\psi_r: X_r \rightarrow X_{r-1}$.
The information of which fibres of 
$\psi_r$ are isomorphic can be organized functorially
by the functor $\Isom_r$, which is representable 
by a scheme $I_r$ locally of finite type over
$X_{r-1} \times X_{r-1}$. $I_r$ comes equipped with
source and target morphisms $s, t: I_r \rightarrow X_{r-1}$.
In the case $S=\P^2$, $\psi_r: X_r \rightarrow X_{r-1}$ together
with $s$ and $t$ solve the moduli problem for
$\F_r$ in the sense explained above (cf. \cite[II, III]{Har88}).
We shall now extend a little this knowledge by relating
Harbourne's  results to Enriques diagrams.

\begin{Pro}
\label{diagsfam}
Let $T$ be a scheme of finite type over $k$,
and 
$$
\begin{CD}
Y_r \overset{\pi_r}\longrightarrow Y_{r-1} \longrightarrow
  \dots \overset{\pi_1}\longrightarrow Y_0=S \times T.
\end{CD}
$$
an ordered blowing-up of $r$ $T$--points. 
Then for every irreducible component $T_i$ of $T$, there exist an 
open subset $U_i \subset T_i$ and an ordered Enriques
diagram $\D$ such that for every closed point $t\in U_i$, 
the cluster $K(t)$ has Enriques
diagram $\D$. Moreover, for every closed point $t \in T_i$
let $\D_t$ be the Enriques diagram of $K(t)$. Then,
denoting by $P$ and $P_t$ the proximity matrices
of $\D$ and $\D_t$ respectively, 
$P_t^{-1}P$ has no negative entries.
\end{Pro}
\begin{proof}
It is not restrictive to assume $T$ irreducible.
Since $X_{r-1}$ represents the functor $\B_r$,
attached to the ordered blowing-up there is a
morphism $T \rightarrow X_{r-1}$ given by
$t \mapsto K(t)$. As for every Enriques diagram
$\D$, $\Cl(\D)$ is locally closed in $X_{r-1}$,
there is an open subset $U \subset T$ such
that the Enriques diagram of $K(t)$, $t\in U$,
is constant; let it be $\D$. Then clearly
the image of $T$ in $X_{r-1}$ is contained
in $\overline{\Cl(\D)}\subset \Eff(\D)$, and the claim follows
from \ref{noiguals}.
\end{proof}

\begin{Pro}
\label{diagsfam2}
Let $T$ be a scheme of finite type over $k$, and 
$\psi: X \rightarrow T$ a proper and smooth morphism whose
fibres at closed points are blowing-ups of clusters
of $r$ points of $\P^2$. Then for every $t \in T$ there exist 
an \'etale neighborhood $V$ of $t$ and a morphism
$K: V \rightarrow X_{r-1}$
such that $X_V \rightarrow V$ comes from $\psi_r$
by base change $V \rightarrow X_{r-1}$, i.e.
$\psi|_V=K^{*}(\psi_r)$, and
there are an open dense subset $U$ of $V$ and ordered
Enriques diagrams $\D_t$ and $\D_U$ satisfying
\begin{enumerate}
\item The ordered cluster $K(t)$ has Enriques diagram $\D_t$.
\item For every $t'\in U$, the ordered cluster $K(t')$ has 
Enriques diagram $\D_U$.
\item If $P_t$ and $P$ are the proximity matrices of
$\D_t$ and $\D_U$ then $P_t^{-1}P$ has no negative entries.
\end{enumerate}
\end{Pro}
\begin{proof}
The first part of the claim is an immediate
consequence of the fact that 
$\psi_r: X_r \rightarrow X_{r-1}$ together
with $s$ and $t$ solve the moduli problem for
$\F_r$ as explained above. It is not restrictive
to assume that the \'etale neighborhood $V$ is
irreducible; then, applying that for each 
ordered Enriques diagram $\D$,
$\Cl(\D)$ is locally closed in $X_{r-1}$, one obtains
the existence of an open set $U \subset V$ such that
the ordered clusters of $K(U) \subset X_{r-1}$
have constant Enriques diagram, which we may call $\D_U$.
Clearly one has 
$K(V) \subset \overline{\Cl(\D_U)} \subset \Eff(\D_U)$
and then \ref{noiguals} completes the proof. 
\end{proof}

It is worth noting that the fact that the varieties 
$\Cl(\D)$ do not form a stratification of $X_{r-1}$ 
(see  \ref{nostrata}) has consequences
in this context. Indeed, although under the
conditions of Proposition \ref{diagsfam} we can prove
that $P_t^{-1}P$ has no negative entries,
it is not necessarily true that $\D \rightsquigarrow \D_t$.

\section{Fixing the position of points}
\label{puntfix}

We consider next other relations between points
of clusters (not proximity relations) which determine 
closed subvarieties of $X_{r-1}$.
Sometimes one may be interested in clusters supported at a single
given point, or whose points lie completely or partially on a given 
curve. If $S=\P^2$, one may want to impose that some points of the
cluster are aligned (not specifying on which straight line) or on
a curve of a given degree; on an arbitrary surface we would ask
the points to lie in some unspecified effective divisor of a given
divisor class. All this types of conditions determine closed 
subvarieties of $X_{r-1}$, and therefore of $\Eff(\D)$, $\Cl(\D)$
and $\overline{\Cl(\D)}$.

Let $K_0 \in X_{s-1}$ be a (fixed) cluster of $s \le r$ points,
and 
$$\iota : \{1, 2, \dots, s\} \rightarrow \{1, 2, \dots, r\}$$
an injection preserving order. We are interested in the subset 
$\Cl(\iota K_0)$ of $X_{r-1}$ containing the clusters $K$ 
with $p_{\iota(i)}(K) = p_i(K_0)$, $i=1, 2, \dots, s$.  

\begin{Pro}
For every $K_0, \iota$ as above, there is a closed subset
$\Eff(\iota K_0)$ of $X_{r-1}$ such that $\Cl(\iota K_0)$ is open in 
$\Eff(\iota K_0)$ and the morphism
$$ \begin{CD}
S_{\Eff(\iota K_0)} @>\pi>> S
\end{CD}$$
factors through $S_{K_0}$.
In particular, $\Cl(\iota K_0)$ is locally closed in $X_{r-1}$.
\end{Pro}
\begin{proof}
By induction on $s$.
For $s=0$ there is nothing to prove. For $s\ge 1$,
let $\breve K_0$ be the cluster obtained by dropping the last point $p_r$ and
$t= \iota(s)$. Since $\iota$ preserves order, $t-1 \geq s-1$, so by the induction
hypothesis there is a closed subset $\Eff(\iota \breve K_0) \subset X_{t-2}$
such that the set $\Cl(\iota \breve K_0)$ is open in 
$\Eff(\iota \breve K_0)$ and the morphism
$$ \begin{CD}
S_{\Eff(\iota \breve K_0)} @>\pi>> S
\end{CD}$$
factors through 
$$ \begin{CD}
S_{\Eff(\iota \breve K_0)} @>\pi_0>> S_{\breve K_0} \ .
\end{CD}$$
Let $Y= \pi_0^{-1} (p_s(K))$. The strict transform $X$ in $X_{t-1}$ of
$Y \times_{X_{t-3}} \Eff(\iota \breve K_0)$ is canonically isomorphic
to the blowing-up family $\BF(\Eff(\iota \breve K_0), Y, X_{t-3})$, and 
due to the universal property of blowing up $p_s(K)$ in $S_{\breve K_0}$
the morphism 
$$ \begin{CD}
X @>\pi_{t-1}>> \pi_i^{-1}(\Eff(\iota \breve K_0)) 
   @>\pi_{0}>> S_{\breve K_0}
\end{CD}$$
factors through $S_{K_0}$. It is easy to see that 
$\psi_{r-1,t-1}^{-1}(Y)$ fulfills the claimed conditions.
\end{proof}

\subsection{The $r$-th neighbourhood of a point}
\label{sec:entorn}

One especially interesting set of clusters with 
fixed points and proximities is that containing
all unibranched clusters of $r$ points whose only
root is a given point $p$. We call it 
$Y_{r-1}(p)\subset X_{r-1}$ or
simply $Y_r$. There is an obvious bijection
between $Y_r(p)$ and the set of all points in the
$r$-th neighbourhood of $p$, given by mapping
$K$ to its last point $p_r(K)$. Recall that for
unibranched clusters there is just one admissible
ordering, so there is no need to distinguish 
between ordered and unordered unibranched clusters.
We are now going
to see that $Y_{r-1}(p)$ is a smooth projective 
variety, and it is not difficult to see that given two
points $p$, $p'$, even from different surfaces,
the varieties $Y_{r-1}(p)$ and $Y_{r-1}(p')$
are isomorphic. So we may call $Y_{r-1}$ simply the
$r$-th neighbourhood of a point on a smooth
surface. These varieties and their immersions in 
$\P^N$ were a classical object of study
(see \cite{Sem54}). For example J. G. Semple
proved that the variety of points of the
second neighbourhood of a point (i.e., $Y_1$) is
a surface isomorphic to a quintic scroll in
$\P^6$. The geometry of these varieties is 
relevant in the analytic classification of
curve singularities (see \cite{Roe99}).
In chapter \cite{Roe?2} a closer study was made
of the subvarieties of $Y_r$ corresponding to some
Enriques diagrams and their relation to the Hilbert
scheme of points of $S$.

\begin{figure}
  \begin{center}
    \mbox{\includegraphics{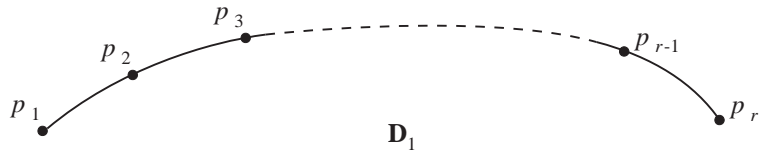}}
    \caption{Enriques diagram $\D_1$.}
    \label{D1}
  \end{center}
\end{figure}

\begin{Pro}
\label{smoothent}
  For every $r \ge 0$, $Y_r$ is a closed smooth 
rational irreducible subvariety
of $X_r$, of dimension $r$. 
\end{Pro}
\begin{proof}
  By its own definition, $Y_r$ can be described
 using the notations introduced
above, $Y_r=\Eff(\D_1(r)) \cap \Cl(\iota p)$,
where $\D_1(r)$ is the Enriques diagram shown in figure 
\ref{D1} and $\iota(1)=1$. The proof goes by
induction on $r$. In the case $r=0$ we have
$Y_0=\{p\}$ and all claims are obvious. For
$r>0$, we give a recursive construction for $Y_r$ and 
$S_{Y_r}=(X_{r+1})_{Y_r}$. We have $Y_{-1}=\Spec{k}$,
$Y_0=\{p\}$ and $S_{Y_{-1}}=S$, whereas for $r \ge 0$,
$$
S_{Y_r}=\BF(S_{Y_{r-1}}, Y_r, Y_{r-1}) \, ,
$$
and $Y_{r+1}$ is the exceptional divisor of the
blowing-up family (because of \ref{construccioXD}).
Now all claims follow from the induction hypothesis.
\end{proof}

\subsection{Points on given curves}

Let now $D$ be a divisor on $S$ and $\L \subset |D|$ a
linear system of curves on $S$. Let $\m=(m_1, m_2, \dots, m_r)$
be a system of multiplicities. We are interested in the subset
$\Eff(\L, \m)$ of $X_{r-1}$ containing 
the clusters $K$ such that there is a curve $C \in \L$ 
going through $(K,\m)$.

\begin{Pro}
\label{clustersistemalineal}
$\Eff(\L,\m)$ is closed in $X_{r-1}$.
\end{Pro}

\begin{proof}
From \cite[4]{KP99} we know that the incidence variety 
$$
V=\left\{ (K, C) \in X_{r-1}\times \L \ | \ C 
\text{ goes through }(K,\m) \right\}
$$
is closed in $X_{r-1}\times \L$.
Then the claim is immediate since $\Eff(\L,\m)$ is
just the projection of $V$ on the first factor.
\end{proof}

Note that the result of \cite{KP99} is much more general,
because it deals with algebraic families of curves
on smooth families of surfaces; we use only a particular
case of it, namely when the surface is projective and fixed
and the family of curves is linear.

\smallskip

Finally, we consider the set of all
clusters in $X_{r-1}$ which have their points (or some
of them) on a given curve $C \subset S$.
Let $C \subset S$ be a curve, and
$I \subset \{1, 2, \dots, r\}$ a set of indices,
fixed for the rest of the chapter.
Let $\Cl(C,I)$ be the set of all clusters $K$ 
such that $p_i(K)$ belongs to $C$ (or its strict 
transform) for all $i$ in $I$.

\begin{Pro}
\label{sobrecconstr}
$\Cl(C,I)$ is constructible.
\end{Pro}

To prove Proposition \ref{sobrecconstr}, we need
the following lemma. Let $\L$ be the linear system 
that consists of the single curve $C$ and, for every system of
multiplicities $\m$, let $\Cl(C,\m) \subset \Eff(\L,\m)$
be the set of all clusters $K$ such that
$\mult_{p_i(K)}(C)=m_i$ for all $i$. 

\begin{Lem}
\label{loccclo}
$\Cl(C,\m)$ is open in $\Eff(\L,\m)$. In
particular, $\Cl(C,\m)$ is locally closed in $X_{r-1}$.
\end{Lem}
\begin{proof}
Let $M \subset \Z^r$ be the following set of
systems of multiplicities:
$$
M=\left\{ \m' \ne \m \, | \, \Cl(C,\m') \cap \Eff(\L,\m) \ne \emptyset \right\}.
$$
$M$ is finite, because the multiplicity of the points
of $C$ is bounded, so
the systems of multiplicities $\m'$ such that
$\Cl(\L,\m') \ne \emptyset$ are a finite set,
and this set contains $M$. We will prove that 
$$
\Eff(\L,\m) \setminus \Cl(C,\m) =
\bigcup_{\m' \in M} (\Eff(\L,\m) \cap \Eff(\L,\m')),
$$
and then the claim will follow since,
after \ref{clustersistemalineal}, $\Eff(\L,\m) \cap \Eff(\L,\m')$
is closed for all $\m$ and $\m'$.

The inclusion 
$$
\Eff(\L,\m) \setminus \Cl(C,\m) \subset
\bigcup_{\m' \in M} (\Eff(\L,\m) \cap \Eff(\L,\m')),
$$
is clear, because if $K \in \Eff(\L,\m) \setminus \Cl(C,\m)$
then defining $\m'=(m_1, m_2, \dots, m_r)$ with 
$m_i=\mult_{p_i(K)}(C)$ we have 
$$
K \in \Eff(\L,\m) \cap \Cl(\L,\m') \subset \Eff(\L,\m) \cap \Eff(\L,\m')
$$
and $\m' \in M$. We shall be done if we prove that,
for every $\m' \in M$,
$$ \Eff(\L,\m) \cap \Eff(\L,\m') \subset
\Eff(\L,\m) \setminus \Cl(C,\m),
$$
or, equivalently, $\Eff(\L,\m') \cap \Cl(C,\m) = \emptyset$
(recall that $\Cl(C,\m) \subset \Eff(\L,\m)$.
So let $\m'=(m_1, m_2, \dots, m_r) \in M$, and let 
$K_1$ be such that $\mult_{p_i(K_1)}(C)=m'_i$ for all $i$
($K_1$ exists by the definition of $M$).
We have to see that for every $K_2 \in \Eff(\L,\m')$,
$K_2 \not \in \Cl(C,\m)$.

Let $j$ be the least index such that 
$m_j \ne m'_j$, and let $S_1=S_{K_1^{j-1}}$
and $S_2=S_{K_2^{j-1}}$
be the surfaces obtained from $S$ blowing
up the first $j-1$ points of $K_1$ and
$K_2$ respectively. As $K_1 \in \Eff(\L,\m)$, clearly
the virtual transform $C_1$ of $C$ in $S_1$
relative to the system of multiplicities $\m$ 
(which coincides with the strict transform
of $C$, as $m_i=m'_i=\mult_{p_i(K_1)}(C)$ for $i<j$)
has at $p_j(K_1)$ multiplicity 
$m'_j=\mult_{p_j(K_1)}(C) \ge m_j$, and we know
$m'_j \ne m_j$. Therefore $m'_j > m_j$.

Now $K_2 \in \Eff(\L, \m')$, so 
the virtual transform $C_2$ of $C$ in $S_2$
relative to the system of multiplicities $\m$ 
(which coincides with the virtual transform
relative to $\m'$, as $m_i=m'_i$ for $i<j$)
has at $p_j(K_2)$ multiplicity at least
$m'_j >m_j$, i.e. $\mult_{p_j(K_2)}(C) \ne m_j$, 
therefore $K_2 \not\in \Cl(C,\m)$, as wanted.
\end{proof}

\begin{proof}[Proof of \ref{sobrecconstr}]
Let now $M \subset \Z^r$ be the following set of
systems of multiplicities:
$$
M=\{ (m_1, m_2, \dots, m_r) \in \Z^r \, | \, m_i >0 \, \forall i\in I\}
$$
The set $\Cl(C,I)$ can be described as
$$\Cl(C,I)= \bigcup_{\m \in M} \Cl(C,\m),$$
and only a finite number of the $\Cl(C,\m)$ are nonempty.
Therefore $\Cl(C,I)$ is a union of sets which after \ref{loccclo}
are locally closed, so $\Cl(C,I)$ is constructible.
\end{proof}

\begin{Cor}
\label{sobrecorba}
Assume $C$ is smooth, and let $\Cl_0(C,I)\subset \Cl(C,I)$ 
be the set of all clusters $K$ such that
$p_i(K)$ belongs to $C$ if and only if $i$ belongs to $I$.
Then $\Cl_0(C,I)$ is locally closed in $X_{r-1}$.
\end{Cor}

\begin{proof}
  As $C$ is smooth, $p_i(K)$ belongs to $C$
if and only if $\mult_{p_i(K)}(C)=1$. Therefore
$\Cl_0(C,I)=\Cl(C,\m)$, where 
$\m=(m_1, m_2, \dots, m_r)$ is the system 
of multiplicities defined by
$$
m_i = 
\begin{cases}
  1 & \text{ if }i\in I, \\
  0 & \text{ if }i\not\in I.
\end{cases}
$$
Now the claim is immediate, after lemma \ref{loccclo}.
\end{proof}

\bibliographystyle{amsplain}
\bibliography{Biblio}
\end{document}